\documentclass[trsc,nonblindrev]{informs3}
\OneAndAHalfSpacedXI 
\usepackage{natbib}
 \bibpunct[, ]{(}{)}{,}{a}{}{,}%
 \def\bibfont{\small}%

\usepackage{algorithm,algpseudocode}
\usepackage{algcompatible}
\usepackage{subcaption}
\usepackage{caption}
\usepackage{graphicx,amsmath}
\usepackage{comment}
\usepackage{mathtools}
\usepackage{booktabs}
\usepackage{multirow}
\usepackage{graphicx}
\usepackage{float}
\newfloat{algorithm}{t}{lop}
\usepackage{tikz}
\usepackage{pgfplots} 
\usepackage{pgf}
\usepackage{pgfplotstable}
\usepackage{lmodern}
\usepackage{import}
\tikzset{every picture/.style={line width=0.75pt}}
\usepackage{graphicx}
\usepackage{bbm}
\usepackage{tablefootnote} 
\usepackage{csquotes}
\pgfplotsset{compat=1.18} 

\TheoremsNumberedThrough  
\EquationsNumberedThrough   
\MANUSCRIPTNO{none} 

\begin{document}
\RUNAUTHOR{Celik et al.}

\RUNTITLE{Two-Step Benders Decomposition}

\TITLE{Exact Two-Step Benders Decomposition for the Time Window Assignment Traveling Salesperson Problem} %

\ARTICLEAUTHORS{
\AUTHOR{\c{S}ifa \c{C}elik$^{a}$, Layla Martin$^{a,b}$, Albert H. Schrotenboer$^{a,b}$, Tom Van Woensel$^{a,b}$}
\AFF{$^a$School of Industrial Engineering, Eindhoven University of Technology, 5612AZ Eindhoven, The Netherlands}
\AFF{$^b$Eindhoven AI Systems Institute, Eindhoven University of Technology, 5612AZ Eindhoven, The Netherlands} \AFF{\EMAIL{s.celik@tue.nl}, \EMAIL{l.martin@tue.nl}, \EMAIL{a.h.schrotenboer@tue.nl}, \EMAIL{t.v.woensel@tue.nl}}
} 

\ABSTRACT{ 
Next-day delivery logistics services are redefining the industry by increasingly focusing on customer service. Each logistics service provider's challenge is jointly optimizing time window assignment and vehicle routing for such next-day delivery services. To do so in a cost-efficient and customer-centric fashion, real-life uncertainty, such as stochastic travel times, needs to be incorporated into the optimization process. This paper focuses on the canonical optimization problem within this context: the Time Window Assignment Traveling Salesperson Problem with Stochastic Travel Times (TWATSP-ST). It belongs to two-stage stochastic mixed-integer programming problems with continuous recourse. We introduce Two-Step Benders Decomposition with Scenario Clustering (TBDS) as {an} exact solution methodology for solving such stochastic programs. The method utilizes a new two-step decomposition along the binary and continuous first-stage decisions and introduces a new scenario-retention strategy that combines and generalizes state-of-the-art Benders approaches and scenario-clustering techniques. Extensive experiments show that TBDS is superior to state-of-the-art approaches in the literature. It solves TWATSP-ST instances with up to 25 customers to optimality. It provides better lower and upper bounds that lead to faster convergence than existing state-of-the-art methods. We use TBDS to analyze the structure of the optimal solutions. By increasing routing costs only slightly, customer service can be improved tremendously, driven by smartly alternating between high- and low-variance travel arcs to reduce the impact of delay propagation throughout the executed vehicle route.}

\HISTORY{}

\maketitle

\section{Introduction}\label{Sec:Introduction}
Logistics service providers (LSPs) face the daily challenge of delivering goods to customers cost-effectively. For years, their emphasis has been on minimizing operational costs. However, competition and globalization require shifting from a purely cost-minimizing mindset to balancing cost and service { quality}. One approach to increase customer service is to assign and communicate a delivery time window well in advance; for example, the evening before delivery in case of next-day parcel delivery. For LSPs, this implies that vehicle routes and time windows must be decided upon simultaneously. Providing a reliable time window for customers is essential, as they will organize their dependent planning processes based upon these in-advance communicated time windows (\cite{spliet2015time, jabali2015self}). For instance, the workforce may be scheduled in case of retail replenishment, or production plans may need to be altered in case of repairs by technicians.

In practice, providing reliable time window assignment and vehicle routing decisions is complex as travel times are inherently stochastic. Consequently, the vehicle routing and time window assignment cannot be considered independently; a slightly longer vehicle route could result in a much better time window assignment performance, and vice versa. The canonical optimization problem addressing this joint issue is the Time Window Assignment Traveling Salesperson Problem with Stochastic Travel Times (TWATSP-ST), which asks for jointly minimizing vehicle routing and time window assignment costs considering stochastic travel times. The TWATSP-ST aligns with the active research stream on vehicle routing that considers time window assignment as an integral part of the optimization problem (see, e.g., \cite{jabali2015self}; \cite{spliet2018time}; \cite{vareias2019assessing}), rather than solely adhering to exogenously given time window constraints (see, e.g., \cite{paradiso2020exact}; \cite{wolck2022branch}). 

The TWATSP-ST { can be formulated as} a two-stage mixed integer stochastic program. Two-stage stochastic programs involve initial first-stage decisions before uncertainty is resolved and second-stage recourse decisions after uncertainty is realized. The objective is to minimize the expected cost associated with both sets of decisions. {Regarding the TWATSP-ST, its two-stage stochastic program comprises} binary and continuous first-stage decision variables (i.e., routing and time window assignment) and continuous second-stage decision variables related to time window violation. {Currently, no exact solution method exists for this problem, and we, therefore, lack an understanding of how to combine time window assignment and vehicle routing optimally.}
{ In this paper, we present a new exact solution methodology for solving the TWATSP-ST, namely \textit{Two-Step Benders Decomposition with Scenario Clustering} \textit{(TBDS)}. Our TBDS is an exact solution method specifically designed to solve the TWATSP-ST. Extensive computational experiments demonstrate that TBDS outperforms state-of-the-art approaches like Benders dual decomposition and partial Benders decomposition. We can optimally solve instances with up to 25 customers and show how integrating the two decisions (compared to previous sequential heuristics) improves the decision quality comprised of travel costs, time window width, and time window violations. Our approach also automatically caters to delay propagation throughout the execution of a vehicle route by solving the TWATSP-ST to optimality. This has either not been accounted for \citep{vareias2019assessing}, or has only been studied with exogenously given time windows \citep{ehmke2015ensuring}. Finally, we numerically show that the structure of the optimal solution smartly divides the high-variance travel-time arcs throughout the route.}

{ Our new TBDS method explicitly exploits the fact that we have both continuous and binary first-stage decision variables. It is general and can be applied to other two-stage stochastic programs involving binary and continuous first-stage decisions. Benders decomposition-inspired approaches, going back to the seminal works of \cite{Benders1962partitioning}, \cite{danzig}, and \cite{van1969shaped}, are commonly utilized to decompose the problem across scenarios. Recent significant advancements, including Benders Dual Decomposition \citep{BDD}, tailored Lagrangian cuts \citep{chen2022generating}, and Partial Benders Decomposition \citep{PartialBenders}, renewed interest in these techniques. 

The effectiveness of our novel TBDS stems from introducing two new, fundamental ideas that each independently already significantly further current Benders methods:} Firstly, we introduce a new two-step decomposition approach that handles the first-stage binary and continuous decision variables to generate optimality and feasibility cuts. Secondly, TBDS proposes a new scenario-retention strategy that integrates recent advancements in scenario-clustering techniques for stochastic programming \citep[see, e.g.,][]{keutchayan2020scenario,clustering}, hereby generalizing the principles of partial Benders decomposition.

The first key concept of TBDS involves decomposing the two-stage stochastic program into a master problem and $N + 1$ subproblems, where $N$ represents the number of scenarios. The initial subproblem corresponds to the first step of the two-step decomposition, focusing solely on the binary decision variables from the master problem and considering a single continuous subproblem. This first step offers a notable advantage by producing a more robust solution for the continuous first-stage decisions, which is then utilized in combination with the binary first-stage solution in the subsequent~$N$ single-scenario subproblems. Consequently, this two-step approach generates significantly stronger multi-optimality cuts, leading to faster convergence to the optimal solution. We theoretically show the benefit of this two-step decomposition.

The second key concept of TBDS is a new scenario-retention strategy that generalizes partial Benders decomposition by incorporating \textit{representative} scenarios into the master problem. These representative scenarios are carefully selected, allowing us to optimize our decisions in the first stage before the uncertainty is observed. The selection of the scenario set involves a trade-off: a larger set of scenarios reflects underlying uncertainty better but increases computational complexity.
There is growing interest in scenario generation and reduction methods to address this trade-off. These methods can be categorized as either distribution-driven, such as those proposed by, e.g., \cite{kleywegt2002sample,henrion2009scenario,pflug2015dynamic}, or problem-driven, as discussed by, e.g., \cite{henrion2018problem}, \cite{clustering}, and this paper. 

The remainder of the paper is organized as follows. In Section \ref{Sec:Problem_Statement0}, we introduce the TWATSP-ST and provide the mathematical model. Next, Section \ref{Sec:TBDS} presents the TBDS methodology, explaining the two key innovations within TBDS. Section \ref{Sec:Problem_Statement} specifies how TBDS will be used to solve the TWATSP-ST. Section \ref{Sec:Comp_Results} provides computational results on the performance of TBDS and managerial results associated with solving the TWATSP-ST. We conclude this paper and provide avenues for future research in Section \ref{Sec:Conclusion}. 

\section{Problem Definition} 
\label{Sec:Problem_Statement0}

The Time Window Assignment Traveling Salesperson Problem with Stochastic Travel Times (TWATSP-ST) concerns the a-priori joint optimization of a vehicle route and the assignment of time windows to a given set of customers in the presence of travel time uncertainty. 

The TWATSP-ST is defined on a graph $G = \left(V,A\right)$, where $V=\{0,\dots, n\}$ is the set of nodes and $A := \{\left(i,j\right)\in V \times V: i\neq j\}$ is the set of arcs. Node $0$ acts as the depot at which the vehicle starts its tour, and all other nodes represent customers. The vehicle has a shift duration of length $T$. Each arc $(i,j) \in A$ has a known distance $d_{ij}\geq 0$. 
Each customer $i\in V^+ := V\setminus \{0\}$ faces a deterministic service time $s_i \geq 0$. We assume travel times over the arcs are stochastic with known distribution. Let \textit{$\xi= \{t_{ij}\}_{(i,j)\in A}$} represent the stochastic travel time vector on a scenario sample space $\Omega$. $\xi_\omega = \{t_{ij\omega}\}_{(i,j)\in A}$ denotes the particular realization of the travel time over arc $(i,j) \in A$ in a scenario $ \omega \in \Omega$. In this way, we intrinsically cater for delay propagation, unlike \cite{vareias2019assessing}. 

In this context, the TWATSP-ST makes two interdependent, a priori decisions: i) A vehicle route visiting all customers in $V^+$, starting and ending at the depot and ii) {a time window assignment { $y_i= \{y^s_i,y^e_i\}$ }  for each customer $i\in V^{+}$}. The objective of the TWATSP-ST is to minimize a weighted sum of expected earliness and lateness at the assigned time windows, the width of the assigned time window, the expected shift overtime of the vehicle, and the total distance the vehicle travels.

{We encode the routing decision with variables $x_{ij} \in \{0, 1\}$ for all $(i,j) \in A$. Together with the time window assignment variables $y$, these form the first-stage decisions in the TWATSP-ST. After the realization of the uncertainty, we can determine the earliness $e_{i\omega}$ and lateness $l_{i\omega}$ to the arrival time of customer $i$ relative to the assigned time window $y_i$ for each customer $i \in V^+$, the departure time $w_{i\omega}$ from a customer that is the time after customer $i$ is served, and the shift overtime $o_{\omega}$ relative to the shift duration length $T$. We define $z_{\omega} = \{w_{i\omega}, e_{i\omega}, l_{i\omega}, o_{\omega}\}_{i\in V^{+}}$ to encapsulate { all} second-stage decision variables. 
}
 
Then, we formulate the TWATSP-ST as the following two-stage stochastic mixed-integer program with continuous recourse. 
\begin{align}
\min \quad & \sum_{i \in V}{\sum_{j \in V\setminus \{i\}} {d_{ij}x_{ij}}} + \sum_{j \in V^+}{\sigma \left(y_{i}^e-y_{i}^s\right)} + \sum_{\omega \in \Omega}p_\omega Q(x,y,\omega) \label{TWATSP:obj_function} \\
\text{s.t.} \quad & \sum_{i \in V\setminus\{j\}} {x_{ij}} = \sum_{i \in V\setminus\{j\}} {x_{ji}} = 1 & \forall j \in V, \label{TWATSP:cons_flow}\\
&\sum_{i\in S}\sum_{j \notin S }{x_{ij}} \geq 1 & S \subseteq V, 1\leq |S| \leq |V^+|, \label{TWATSP:cons_subtour_elimination}\\
& y_{i}^e-y_{i}^s \geq s_i & \forall i \in V^+, \label{TWATSP:cons_timewindow_limit} \\
&x_{ij} \in \{0,1\}  &\forall (i,j) \in  A, \label{TWATSP:domain_ex}\\
&y_{i}^e, y_{i}^s \in \mathbb{R_+} & \forall i \in V^+. \label{TWATSP:domain_tw}
\end{align}
Here, $\sigma \geq 0$ is an exogenously set weight factor. Objective \eqref{TWATSP:obj_function} consists of the travel and weighted time window assignment costs and the second-stage recourse function $Q(x,y,\omega)$. Constraints \eqref{TWATSP:cons_flow} ensure that the vehicle visits each customer exactly once. Constraints \eqref{TWATSP:cons_subtour_elimination} eliminate sub-tours, and constraints \eqref{TWATSP:cons_timewindow_limit} ensure that the time window assignment respects the service time at each customer.

The recourse function $Q(x,y,\omega)$, as part of the Objective  \eqref{TWATSP:obj_function}, equals the expected shift overtime cost and the expected incurred earliness and lateness cost due to the time window assignment and vehicle routing decision.
{
\begin{align}
Q(x, y,\omega):= \min \quad & {\sum_{j \in V^+}\phi\left( e_{j_\omega}+l_{j\omega} \right) + \psi o_\omega} \span \label{TWATSP:obj_second_stage}\\
\text{s.t. } \quad & w_{j\omega} \geq  w_{i\omega} + t_{ij\omega}+  s_j  - \left(1-x_{ij}\right)M & \forall  (i,j ) \in A,  \label{TWATSP:cons_departure}\\
& e_{j\omega} \geq y_{j}^s - w_{j\omega}- s_{j} & \forall j \in V^+, \label{TWATSP:cons_earliness}\\
&l_{j\omega} \geq w_{j\omega} - y_{j}^e & \forall j \in V^+, \label{TWATSP:cons_lateness}\\
&o_\omega \geq w_{j\omega} +t_{j0\omega} -T  & \forall  j \in V^+, \label{TWATSP:cons_overtime}\\
&w_{0\omega} = t_0, &   \label{TWATSP:cons_departure_0}\\
&w_{j\omega},e_{j\omega},l_{j\omega} \in \mathbb{R_+} & \forall j \in V^+, \label{TWATSP:domain_e_l_w}\\
&o_\omega \in \mathbb{R_+}. & \label{TWATSP:domain_o}
\end{align}
}
Here, $ \psi$ and $\phi$ are {exogenous} weight factors for the earliness and lateness concerning the time window assignment and the overtime of the vehicle, respectively.
The objective function of the second-stage problem \eqref{TWATSP:obj_second_stage} is a function of first-stage variables ${(x,y)}$ and a scenario $\omega \in \Omega$.
Constraints \eqref{TWATSP:cons_departure}~-~\eqref{TWATSP:cons_overtime} determine each customer's departure time, earliness, lateness, and overtime. We set $M$ in \eqref{TWATSP:cons_departure} - \eqref{TWATSP:cons_lateness} equal to the longest total travel time among all scenarios. We require tours to start from the depot at a predetermined time $t_0$ (Constraints \eqref{TWATSP:cons_departure_0}). Constraints \eqref{TWATSP:domain_ex}, \eqref{TWATSP:domain_tw}, \eqref{TWATSP:domain_e_l_w} and \eqref{TWATSP:domain_o} define the variable domain. Model \eqref{TWATSP:obj_function} - \eqref{TWATSP:domain_o} is a two-stage stochastic mixed-integer program with continuous recourse that fully defines the TWATSP-ST.


\section{Exact Two-Step Benders Decomposition with Scenario Clustering (TBDS)} \label{Sec:TBDS}

{This section introduces our TBDS solution method. Before we explain the main novelties of TBDS in Sections 3.1 and 3.2, we first briefly introduce Benders decomposition, outlining the essential formulations that TBDS builds upon. To ensure a clear understanding and compact presentation of TBDS for solving the TWATSP-ST, we present it using generic scenario-based two-stage stochastic programming notation. Scenario independence (reflecting the realization of random vectors) ensures that a two-stage stochastic program is decomposable over scenarios. The result is a single master problem and one or multiple subproblems depending on the selected decomposition among the scenarios. These subproblems represent the second-stage recourse problem of the decomposed scenarios $\omega \in \Omega$. Here, $\Omega$ is a finite set of scenarios or realized random events. Let {{$\xi$}} represent a random vector on a finite scenario sample space $\Omega$. Specifically, $\xi_\omega$ denotes the particular realization of the random vector in scenario $\omega \in \Omega$, with $p_\omega$ as the associated probability.}

{
In this paper, we consider general two-stage stochastic mixed-integer programs with relatively complete continuous recourse:
\begin{equation} \label{model:MIP}
 \min_{x,y,z} \{   c^\intercal x + d^\intercal y + \sum_{\omega \in \Omega}{p_{\omega}}f_{\omega}^\intercal  z_\omega \mid
W_\omega x + T_\omega y + S_\omega z_\omega \geq h_\omega \ \forall \ \omega \in \Omega,
{ (x, y )\in \mathcal{X}}, 
z_\omega \in \mathbb{R}_+^m  \ \forall \  \omega \in \Omega.\},
\end{equation}    

where $x \in \mathbb Z_{+}^{n_1}$ denotes a vector of first-stage integer decision variables,  $y \in \mathbb R_{+}^{n_2}$ denotes a vector of first-stage continuous decision variables, and $z_{\omega} \in \mathbb{R}^m$ denotes a vector of continuous second-stage decision variables for each scenario $\omega \in \Omega$. Furthermore, $c \in \mathbb{R}^{n_1}, d \in\mathbb{R}^{n_2},$ and $f_\omega \in\mathbb{R}^m$ are the cost vectors of the first- and second-stage decisions, $\ W_\omega \in \mathbb{R}^{\ell\times n_1 }, T_\omega \in \mathbb{R}^{\ell \times n_2},$ and $ S_\omega \in \mathbb{R}^{\ell \times m}$ are the technology matrices capturing the constraints that link the first- and second-stage decisions, and $h_\omega \in \mathbb{R}^\ell$ is the right-hand side vector associated with the linking constraints.  We denote the first-stage constraints and their domains compactly as { $ (x,y) \in\mathcal{X}~:=~\{ (x,y) \in \mathbb{Z}_+^{n_1}\times\mathbb{R}_+^{n_2} \mid Ax + By \geq a\}$} where $A \in \mathbb{R}^{k\times n_1},\  B \in \mathbb{R}^{k\times n_2},$ and $a \in \mathbb{R}^{k}$. 
We assume that the above problem is primal feasible and bounded. 
}
\begin{remark}
    Comparing the general two-stage stochastic program with the TWATSP-ST, the $x$ variables are the first-stage routing decisions, and the $y$ variables are the first-stage time window assignment decisions. The $z_\omega$ variables are the second-stage decision variables that capture the time window exceedances at each customer, the vehicle departure time, and the vehicle overtime for each travel time scenario $\omega \in \Omega$. 
\end{remark}

To use Benders decomposition for solving problem \eqref{model:MIP}, we introduce an auxiliary decision variable $\theta_\omega$ that approximates the second-stage cost for a scenario $\omega \in \Omega$ and subsequently decompose the problem into a master problem and multiple subproblems. We define the master problem (MP) as  
\begin{align}
MP= \min_{x,y,\theta_\omega}  \quad & c^\intercal x + d^\intercal y  + \sum_{\omega \in \Omega}{p_{\omega}}\theta_{\omega}&  \label{MP:obj_function}\\
\text{s.t.} \quad & \theta_\omega \geq \bar{\theta}_\omega & \forall \ \omega \in \Omega, \label{MP:cons_theta_lower_bound}\\
& (x,y) \in \mathcal{X},  \label{MP:cons_domain_xy}  \\
& \theta_\omega \in \mathbb{R} &  \forall \ \omega \in \Omega, \label{MP:cons_domain_theta_omega} \\
& \text{Feasibility and Optimality Cuts}.   \label{focuts}
\end{align}

where $ \bar{\theta}_\omega$ is a lower bound on $\theta_\omega$ to avoid unboundedness, for each $\omega \in \Omega$. The Feasibility and Optimality Cuts \eqref{focuts} bound the value of $\theta_\omega$ in terms of linear equations in $(x,y)$ by the logic of Benders decomposition. These cuts are generated in the following manner.

The solutions of the MP define lower bounds for the true two-stage stochastic program \eqref{model:MIP}. We will solve the MP in a branch-and-cut manner. That is, at each integer branch-and-bound node, the MP solution values $x^*$ and $y^*$ are fixed in the second-stage subproblem $SP(x^*,y^*,\omega)$ for each scenario~$\omega$:
\begin{align} \label{model:SP}
SP(x^*,y^*) := &\min_{x_\omega,y_\omega,z_\omega} \left\{f_{\omega}^\intercal  z_\omega \mid  W_\omega x_\omega + T_\omega y_\omega + S_\omega z_\omega \geq h_\omega, \right. \\ &\hspace{1.4cm} \left. x_\omega = x^*, y_\omega = y^*,\ x_\omega \in \mathbb{R}^{n_1}_+, y_\omega \in \mathbb{R}^{n_2}_+,  z_\omega \in \mathbb{R}^m_+\right \}. \nonumber
\end{align}
Note that the above problem is equivalent to the true two-stage stochastic program \eqref{model:MIP} without the first-stage constraints, as they are redundant due to the constraints $\{ x_\omega = x^*,\  y_\omega = y^*\}$. We directly decompose this subproblem over each scenario and write $x_\omega$ and $y_\omega$ in the constraint set. Benders decomposition then takes the extreme points or the direction of the extreme rays (in case of infeasibility) of the feasible region of the dual problem of \eqref{model:SP} to generate optimality and feasibility cuts, respectively. These cuts will be added to MP and reduce the gap between MP and the true two-stage stochastic program \eqref{model:MIP}. This procedure is followed at each integer branch-and-bound node while solving the MP, ensuring it converges to the optimal integer solution. 

We now make the optimality and feasibility cuts specific. 
For feasible  MP solutions $x^*$ and $y^*$, let $\nu_\omega$ and $\eta_\omega$ denote the dual multipliers to the constraints $x_\omega = x^*$ and $y_\omega = y^*$ in $SP(x^*, y^*, \omega)$, respectively. By solving the dual of the problem \eqref{model:SP}, we retrieve the extreme point $\grave \nu_\omega$ and $\grave \eta_\omega$ and accordingly, a \textit{standard multi-optimality cut} can be derived
\begin{equation} \label{optimalitycut:multicut}
    \theta_\omega \geq  f_{\omega}^\intercal  \grave{z}{_\omega} + (x-\grave{x}_\omega)^\intercal\grave \nu_\omega + (y-\grave{y}_\omega)^\intercal \grave \eta_\omega,
\end{equation} 

where $(\grave{x}_\omega, \grave{y}_\omega, \grave{z}{_\omega})$ is the optimal primal solution of $SP(x^*, y^*, \omega)$.

If problem \eqref{model:SP} is infeasible for the first-stage solutions $x^*$ and $y^*$, indicating that the dual problem is unbounded, we utilize the dual of the problem \eqref{model:SP} and generate feasibility cuts with the unbounded extreme ray of the feasible region of the dual problem to cut the solution $(x^*, y^*)$ off from the MP. To generate the feasibility cut, we solve the following feasibility problem 
\begin{align}
    \min_{x,y,z,\epsilon} \{ \mathbbm{1}^\intercal  \epsilon \mid W_\omega x + T_\omega y + S_\omega z_\omega + \epsilon \geq h_\omega \ \forall \ \omega \in  \Omega, \ x =x^*,  y = y^*,  z_\omega \in \mathbb{R}^m_+, \epsilon \in \mathbb{R}^\ell_+ \}. \label{model:Feas_Problem}
\end{align}

Note that the feasibility problem has not been decomposed over the various scenarios. This generates a feasibility cut of the form 
\begin{equation}
    0 \geq \mathbbm{1}^\intercal \breve{\epsilon} + (x-\breve{x})^\intercal \breve{\lambda} + (y-\breve{y})^\intercal \breve{\beta}, \label{feasibilitycut} 
\end{equation}
{ where $\breve{\epsilon}$, $\breve{x}$ and $\breve{y}$ refer to the values of $\epsilon$, $y$ and $x$ in the optimal solution to the feasibility problem \eqref{model:Feas_Problem}. Here, $\breve{\lambda}$ and $\breve{\beta}$ are dual solutions to the constraints $x = x^*$ and $y = y^*$, respectively and $\mathbbm{1}$ is a vector of ones of size $\ell$. } 

The optimality and feasibility cuts \eqref{optimalitycut:multicut} and \eqref{feasibilitycut} are added to the MP (i.e., equation \eqref{focuts}) if the current MP solution is violating any of those cuts. Otherwise, problem \eqref{model:MIP} is solved to optimality. 
This solution process sums up what is often referred to as the multicut version of the L-shaped method (\cite{birge1988multicut}, see \cite{van1969shaped}, \cite{gendreau1998generalized} for the single cut version wherein the scenario-based cuts generated at each iteration are aggregated into a single cut). 

Adding optimality and feasibility cuts iteratively removes all guidance regarding the continuous first-stage decision variables $y$ concerning the second stage in the first iterations from MP. Regarding the TWATSP-ST, the first-stage time-window assignment solutions $y$ are not guided by the underlying travel time uncertainty. As a result, the MP solutions in early iterations will be poor in terms of their lower approximations to the associated recourse cost encapsulated in the second-stage subproblems. Moreover, these initial MP solutions might also be far from cost-efficient in the second stage. The solution process can become unduly slow since violated cuts are added only after solving the current MP. In the following section, we will discuss two new strategies to enhance the efficiency of the Benders decomposition algorithm for solving two-stage mixed-integer stochastic programs with continuous recourse. First, we will outline a new decomposition strategy to provide a strengthened version of the standard multi-optimality cuts \eqref{optimalitycut:multicut} in Section 3.1. After that, in Section 3.2, we discuss a new scenario retention strategy for improving the MP lower bound.

{
\subsection{Two-Step Decomposition for Strengthened Multi-Optimality Cuts} \label{subsec:TwoStepDecomposition}
{ We now detail our newly proposed two-step decomposition technique to strengthen the standard multi-optimality cut \eqref{optimalitycut:multicut}. To do so, we consider an alternative master problem $MP-2$. 
\begin{align} \label{model:MP2}
\text{MP-2} :=  \min_{x,y, \Theta, \theta} &  c^\intercal  x  + \Theta &\\ \text{s.t. } & \Theta \geq d^\intercal  y  + \sum_{\omega \in \Omega} {p_{\omega}}\theta_{\omega},& \\ & \theta_\omega \geq \bar{\theta}_\omega \ & \forall \ \omega \in \Omega, \\ & (x,y) \in \mathcal{X}, \\ &\text{Feasibility and Optimality cuts,}  \\ & \Theta \in \mathbb{R}_+, \\ &\theta_\omega \in \mathbb{R} \ &\forall \ \omega \in \Omega.
\end{align}
In MP-2, we define the auxiliary decision variable $\Theta$, which is a lower approximation of the cost associated with the continuous first-stage decision variable $y$ and the recourse cost, i.e., $d^\intercal y + \sum_{\omega \in \Omega}{p_{\omega}}f_{\omega}^\intercal  z_\omega $. Clearly, $\sum_{\omega \in \Omega}\theta_\omega$ is a lower approximation of $\sum_{\omega \in \Omega}{p_{\omega}}f_{\omega}^\intercal  z_\omega$, so it follows directly that it is valid to impose $\Theta \geq d^\intercal y + \sum_{\omega \in \Omega}{p_{\omega}}f_{\omega}^\intercal  z_\omega$. Thus, $MP-2$ is equivalent to the original problem~\eqref{model:MIP}.

We define the following \textit{aggregated subproblem} $AP$ to derive optimality cuts.  
\begin{equation}
AP(\bar x) = \min_{x,y,z} \{d^\intercal y + \sum_{\omega \in \Omega} {p_{\omega}}f_{\omega}^\intercal  z_\omega \mid W_\omega x + T_\omega y + S_\omega z_\omega \geq h_\omega \ \forall \ \omega \in \Omega, \ x =\bar x,  y \in \mathcal{Y}(\bar x), z_\omega \in \mathbb{R}^m_+\}, \label{model:AP}
\end{equation}
where $\mathcal{Y}(\bar x) := \{y \in \mathbb{R}^{n_2}_+ | A\bar x + By \geq a \}$ defines the domain of the $y$ variables for AP. Notice that AP considers only the integer master problem solution $\bar x$ as input and optimizes over the continuous first-stage decision variable $y$ with all the second-stage decision variables of the complete scenario set $\Omega$.
We do not suffer from not decomposing AP over each scenario because it is a linear program and, thus, easy to solve computationally. Including all scenarios provides more information, resulting in a relatively ‘good’ first-stage solution of the continuous first-stage decision variable  $y$. }

{ The core idea of the two-step decomposition in TBDS is as follows. Let $\bar y$ be the optimal solution of $AP(\bar x)$ for a feasible solution $\bar x$ of $MP-2$. Then, we consider the subproblem
\begin{align} \label{model:SP2}
SP(\bar x,\bar{y}, \omega) := &\min_{x_\omega,y_\omega,z_\omega} \left\{f_{\omega}^\intercal  z_\omega \mid  W_\omega x_\omega + T_\omega y_\omega + S_\omega z_\omega \geq h_\omega, \right. \\ &\hspace{1.4cm} \left. x_\omega = \bar x, y_\omega = \bar y,\ x_\omega \in \mathbb{R}^{n_1}_+, y_\omega \in \mathbb{R}^{n_2}_+,  z_\omega \in \mathbb{R}^m_+\right \}, \nonumber
\end{align}
which is similar to subproblem  \eqref{model:SP} but takes the optimal solution $\bar{y}$ from $AP(\bar x)$ instead of taking the optimal solution $y^*$ directly from MP-2 (or MP) as an input.
Then consider the optimal MP-2 solution $\bar x$, the optimal solution $\bar y$ to $AP(\bar x)$, and let $\tilde \nu_\omega$ and $\tilde \eta_\omega$ denote the dual multipliers to the constraints $x_\omega = \bar x$ and $y_\omega = \bar y$ in $SP(\bar x, \bar y, \omega)$, respectively. By solving the dual of the problem \eqref{model:SP2}, we retrieve the extreme points $\tilde \nu_\omega$ and $\tilde \eta_\omega$. Accordingly, we derive a  \textit{strengthened multi-optimality cut}:
\begin{equation} \label{optimalitycut:strmulticut}
    \theta_\omega \geq  f_{\omega}^\intercal  {\tilde{z}}{_\omega} + (x-{\tilde{x}}_\omega)^\intercal \tilde \nu_\omega + (y-{\tilde{y}}_\omega)^\intercal \tilde \eta_\omega,
\end{equation}

where $({\tilde{x}}_\omega, {\tilde{y}}_\omega, {\tilde{z}}{_\omega})$ is the optimal primal solution of $SP(\bar x, \bar y, \omega)$.

To summarize the differences between the traditional approach of deriving \textit{multi-optimality cuts} and our two-step approach of deriving \textit{strengthened multi-optimality cuts}, we summarize the two procedures:
\begin{enumerate}
    \item The traditional procedure solves $MP-2$ (or MP) to obtain $x^*$ and $y^*$, then solves $SP(x^*, y^*, \omega)$ to derive a standard multi-optimality cut \eqref{optimalitycut:multicut}. 
    \item Our two-step approach solves MP-2 to obtain $\bar x$, solves $AP(\bar x)$ to obtain $\bar y$, then solves $SP(\bar x, \bar y, \omega)$ to obtain the strengthened multi-optimality cut \eqref{optimalitycut:strmulticut}.
\end{enumerate}

In the following, we provide theoretical insight on why our two-step procedure is beneficial compared to the standard approach. We will characterize the improvement gained by following our two-step decomposition technique to derive the strengthened multi-optimality cut \eqref{optimalitycut:strmulticut} compared to deriving the standard multi-optimality cut \eqref{optimalitycut:multicut}. 
 In the interest of brevity, we use $\mathcal{Z}_\omega := \{z_\omega \in \mathbb{R}^m_+  \mid W_\omega x_\omega + T_\omega y_\omega + S_\omega z_\omega \geq h_\omega \text{ for some } y_\omega \in \mathbb{R}^{n_2}_+, x_\omega \in \mathbb{Z}_+^{n_1}  \}$ and  $\Tilde{\mathcal{X}} = \mathcal{X} \cap\{(x, y) \in \mathbb{Z}_+^{n_1} \times \mathbb{R}^{n_2}_+ \mid W_\omega x + T_\omega y + S_\omega z_\omega \geq h_\omega: \text{ for some }  z_{\omega} \in  \mathbb{R}^m_+ \}$ to represent the feasible region of the $SP(\cdot,\cdot, \omega)$ for $\omega \in \Omega$ and to indicate the domain of the $y$ variables constrained by the $z$  variables, respectively.

The following proposition shows the benefit of taking our two-step decomposition compared to a classic decomposition. That is, we show the improvement gained by restraining the continuous first-stage variables $y$ by taking their solution via $AP(\cdot)$ instead of directly deriving a multi-optimality cut, for arbitrary dual vectors $\hat \nu$ and $\hat \eta$.

\begin{proposition} \label{propososition:stronger_cut}
Let  $ \hat \nu \in \mathbb{R}^{n_1}$ and $\hat \eta \in \mathbb{R}^{n_2}$ be any dual multipliers feasible for $SP(\cdot,\cdot, \omega)$. 
The { strengthened} multi-optimality cut is at least $\Delta \geq 0$ units tighter than the { standard} multi-optimality cut, where
\begin{equation} \label{eq:deltadefinition}
    \Delta =  \min_{\substack{ z_{\omega}  \in conv(\mathcal{Z}) \\  (x_\omega, y_{\omega}) \in conv(\Tilde{\mathcal{X}}) }} \{ f_\omega^\intercal z_\omega   -x_\omega^{\intercal} \hat \nu {- y_\omega}^\intercal  \hat \eta \}  - \min_{\substack{ z_{\omega}  \in conv(\mathcal{Z}) \\  (x_\omega,y_{\omega}) \in conv({\mathcal{X}}) }} \{ f_\omega^\intercal z_\omega   -x_\omega^{\intercal} \hat \nu {- y_\omega}^\intercal \hat \eta \}  .
\end{equation}
\end{proposition}
The intuition of the proposition is that the left minimization term is derived from the strengthened optimality cut that considers $(x_\omega, y_\omega)\in conv(\tilde{\mathcal X})$ which ensures that $y_\omega$ is restricted to be the solution of $AP(\cdot)$ whereas the right minimization term is derived from the standard optimality cut.

Besides using AP \eqref{model:AP} to generate the strengthened multi-optimality cuts, we can also use AP to generate another set of optimality cuts by following the same logic as for deriving cuts based on the problem \eqref{model:AP}. Working in a branch-and-cut scheme, we consider an arbitrary node of the branch-and-bound tree and the corresponding solution $\bar x$ of subproblem AP. Let $\bar \mu $ be the optimal value of the dual multipliers associated with the constraints $x = \bar{x}$ in AP. If AP is infeasible for $\bar x$, we generate the feasibility cut \eqref{feasibilitycut}. On the other hand, if the primal subproblem (AP) returns a feasible solution $(\Bar{y}, \bar{z})$, we derive a so-called \textit{generalized Benders optimality cut}
\begin{equation}
    \Theta \geq d^\intercal \Bar{y} + \sum_{\omega \in \Omega}{p_{\omega}}f_{\omega}^\intercal  \Bar{z}{_\omega} + (x-\Bar{x})^\intercal \bar \mu.\label{optimalitycut:generalized_benders_cut}
\end{equation}
}}
We illustrate and provide insights into the benefits of the proposed cut-generation strategy with a toy example. 
\begin{example} 
Consider the following toy example:
\begin{align*} 
\min_{\substack{x \in \mathbb{Z}_+,z, \\ 2 \leq y \leq 8}} \ \{ y + z:
- 2x -3y + 5z  \geq 17, 
 3y + 2z \geq 10,
2x - z \geq -10,
- 5x + 10y + 2z   \geq 11, 
 x + y + 2z \geq 15\}
\end{align*}
The optimal solution of the problem is $x = 2$, $y = 2$, $z= 5.5$ with a cost of 7.5. In Figure \ref{Fig:Cuts}, we demonstrate the cuts generated in each iteration in the $(x,y,z)$ space over the $ y = 2$ plane. In Figure \ref{Fig:AOC}, we show the generalized Benders optimality cuts, in Figure \ref{Fig:SOC} we show the strengthened multi-optimality cuts, and in Figure \ref{Fig:OC} we show the standard multi-optimality cuts. The color of the cuts refers to the iteration in which they are added, as explained later.

\begin{figure}[!htbp]
     \centering
    \caption{Performance of Different Optimality Cuts in Solution Space of the LP Master Problem}
     \begin{subfigure}{0.32\textwidth}
        \includegraphics[width=\textwidth]{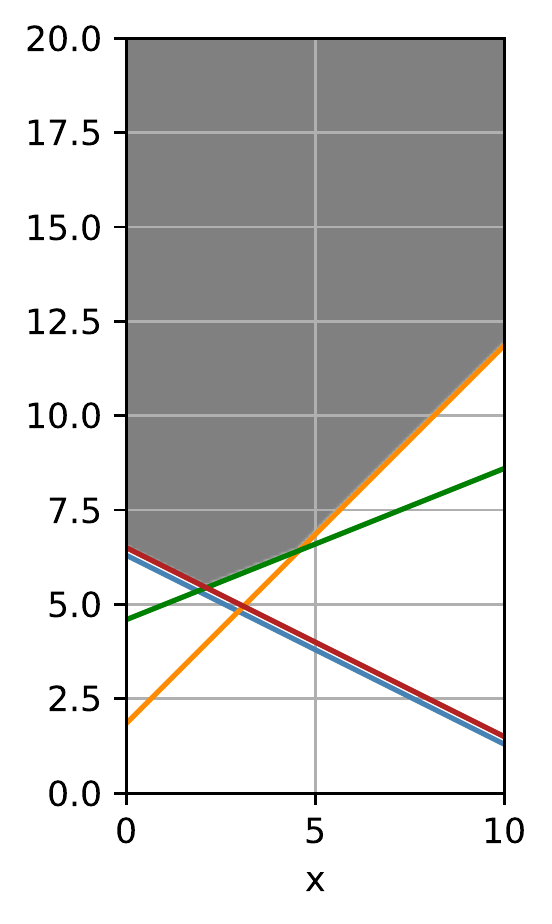}
      \caption{Generalized Benders Optimality Cut}
      \label{Fig:AOC}
     \end{subfigure}
     \hfill
          \begin{subfigure}{0.32\textwidth}
        \includegraphics[width=\textwidth]{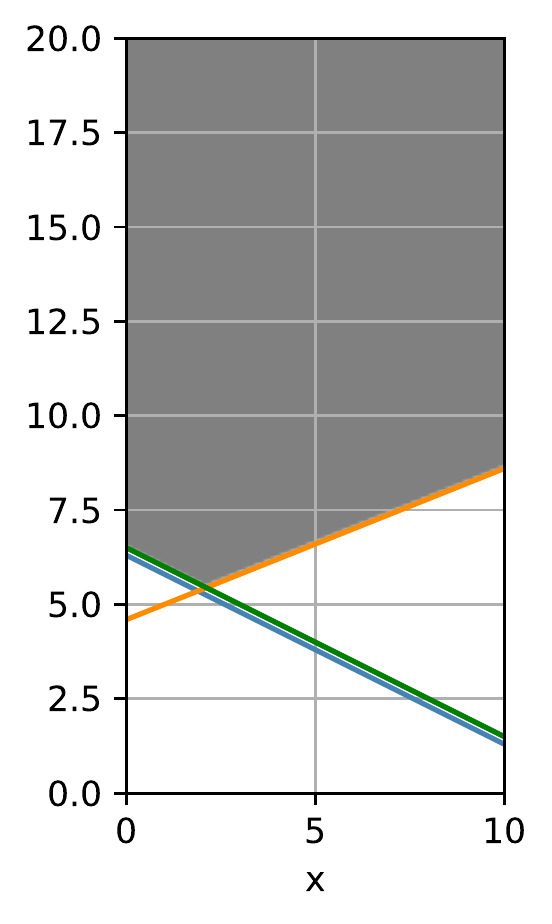}
      \caption{Strengthened Multi-Optimality Cuts}
      \label{Fig:SOC}
     \end{subfigure}
     \hfill
     \begin{subfigure}{0.32\textwidth}
        \includegraphics[width=\textwidth]{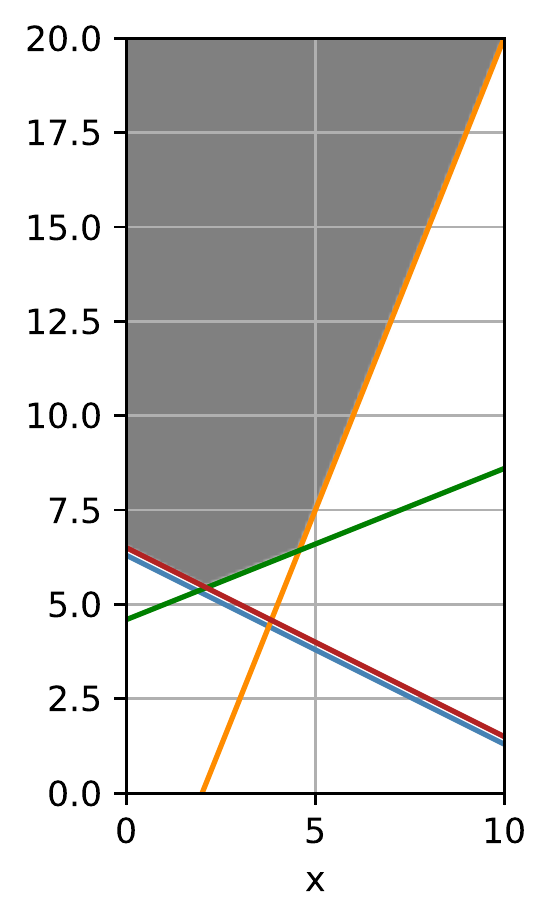}
      \caption{Standard Multi-Optimality Cuts}
      \label{Fig:OC}
     \end{subfigure}
     \label{Fig:Cuts}
\end{figure}

Figure \ref{Fig:OC} shows that relying on standard Benders multi-optimality cuts, we must generate four cuts to reach the optimal solution. As part of our TBDS approach, the strengthened multi-optimality cuts only require three cuts to converge to the optimal solution. The first iteration produces an identical cut for all three optimality cuts, indicated by the blue lines in Figure \ref{Fig:Cuts}. The reason is that the $x$ solution for the first iteration is 0, which causes the $y$ values from the MP and AP to be the same, which is 2. Using the strengthened multi-optimality cuts at the second iteration, we increase the $y$ value from 2 to 4.125, producing the orange cut in Figure \ref{Fig:SOC}). As seen in the figure, this elevates the $x$ value and cuts off a large part of the feasible region. This cut is later implemented in the settings with generalized Benders optimality and standard multi-optimality cuts (see the green lines in Figures \ref{Fig:AOC} and \ref{Fig:OC}). Therefore, they require additional iterations to converge to the optimal value.
\end{example}

\subsection{A New Scenario-Retention Strategy for Improving the MP Lower Bound} \label{subsec:ScenarioSelection} 
The second key innovation within our TBDS method is a new scenario-retention strategy to enhance the MP lower bound. Relatively complete continuous recourse in stochastic mixed-integer programs{, as the TWATSP-ST possesses,} typically entails weak bounds and many (superfluous) iterations of generating cuts as the MP loses all the information with the second-stage variables \citep{rahmaniani2017benders}. We adopt and generalize the idea of partial Benders decomposition \citep{PartialBenders} to overcome this issue. In line with partial Benders decomposition, we include second-stage constraints associated with a subset of scenarios $\Omega_{MP}$, via constraints $\{W_\omega x + T_\omega y + S_\omega z_\omega \geq h_\omega  \forall \ \omega \in \Omega_{MP}\}$ by defining  $ z_\omega \in \mathbb{R}_+^m,  \forall \  \omega \in \Omega_{MP}$ in the MP. However, partial Benders decomposition designs the set $\Omega_{MP}$ using a \textit{row covering strategy} to retain scenarios from $\Omega$ to eliminate many feasibility cuts. Instead, we determine a \enquote*{representative} scenario subset $\Omega_{MP}$. To pick the scenarios from the scenario set \enquote*{in a smart way} to improve the MP lower bound, we propose using scenario clustering as a scenario-retention strategy. We utilize the new promising approach proposed by \cite{clustering}, but in contrast to them, we use it within a Benders decomposition approach. To our knowledge, our paper is the first to introduce this concept for improving the lower bound of the MP within a Benders decomposition context. In the following, we detail the exact scenario-retention strategy.

The representative scenarios in the master problem $\Omega_{MP} = \Omega_{MP}^1\cup \Omega_{MP}^2$ comprise \enquote*{actual} scenarios $\Omega^1_{MP} \subseteq \Omega$ and \enquote*{artificial} scenarios $\Omega^2_{MP}$.  Artificial scenarios are the result of convex combinations of actual scenarios. 

We first detail how to determine $\Omega_{MP}^1$. The first step is to cluster the scenario set $\Omega$. The second step is to select the representative scenarios from each of the clusters, which will together define $\Omega_{MP}^1$. This is formalized in Definition \ref{Def:Clsutering}.

\begin{definition}
     \label{Def:Clsutering}
Let $K$ be the number of clusters. The set $\Omega^1_{MP} = \{r_1, \ldots, r_K$\} is constructed as follows:
\\
\textit{Step 1.} Compute the opportunity-cost matrix $\mathbb{V}= (V_{ij})_{|\Omega| \times |\Omega| } $ where

\begin{equation*}
    V_{ij} = SP((\Hat{x_i},\Hat{y_i}),\omega_j) \qquad \forall (i,j) \in \Omega
\end{equation*}
where $(\Hat{x_i},\Hat{y_i})$ is the optimal solution of the one-scenario subproblem:

\begin{equation*}
    (\Hat{x_i},\Hat{y_i}) \in \argmin_{x,y} SP(x,y,\omega_i) \qquad \forall i \in \Omega
\end{equation*}

\textit{Step 2.} Find a partition of the set $\Omega$ into $K$ clusters $C_1,\ldots,C_K$  and their representative scenarios $r_1 \in C_1,\ldots, r_K \in C_K$ such that 
{
\begin{equation}
    r_{k} \in \argmin \dfrac{|C_k|}{|\Omega|} \left|V_{r_k,r_k} - \frac{1}{|C_k|}\sum_{j \in C_k}V_{r_k,j} \right| \label{eq:ClusteringMiminization}
\end{equation}
}
\end{definition} 

By minimizing the clustering error \eqref{eq:ClusteringMiminization}, we create clusters that best fit the average cost function of all clusters. The reader is referred to \cite{clustering} for more details. 

The above procedure in Definition \ref{Def:Clsutering} adds $K$ representative scenarios into the MP by means of constraints $ \{ W_{\omega} x + T_{{\omega}} y + S_{{\omega}} z_{{\omega}} \geq h_{{\omega}} \forall \ {\omega} \in \Omega^1_{MP} \}$. {In the context of the TWATSP-ST, we observe that this typically leads to clustering travel time realizations that lead to similar routes.}
This reduces the root node optimality gap by improving the linear relaxation of the master problem. However, the value of $K$ should be carefully selected. A $K$ that is too large may lead to overpopulating the master problem, increasing the solution time. 

{ Additionally, we construct artificial scenarios $\Omega_{MP}^2$ based on scenarios in $\Omega_{SP} = \Omega \setminus \Omega_{MP}$. Adding artificial scenarios into the master problem affects the first-stage variables. The goal is to improve the lower bound. To generate artificial scenarios $\Bar{\omega} \in \Omega_{MP}^2$, we use convex combinations of scenarios in $\Omega_{SP}$ and add the constraints $ \{ W_{\Bar{\omega}} x + T_{\Bar{\omega}} y + S_{\Bar{\omega}} z_{\Bar{\omega}} \geq h_{\Bar{\omega}} \forall \ \Bar{\omega} \in \Omega^2_{MP} \}$  by defining  $ z_{\Bar{\omega}} \in \mathbb{R}_+^m,  \forall \  \Bar{\omega} \in \Omega^2_{MP}$ to the MP. 
\begin{definition}
For $ \omega \in \Omega_{SP}$, let $\alpha^{\Bar{\omega}}_{\omega} \geq 0,$ such that $\sum_{\omega \in \Omega_{SP}}{\alpha^{\Bar{\omega}}_{\omega}} =1$. Then, the realization of random vector for artificially generated scenario $\Bar{\omega} \in \Omega_{MP}^2$ is defined as 
\begin{align*}
{W}_{\Bar{\omega}} =\sum_{\omega \in \Omega_{SP}}{\alpha^{\Bar{\omega}}_{\omega}  W_{\omega}}, \quad
{T}_{\Bar{\omega}} =\sum_{\omega \in \Omega_{SP}}{\alpha^{\Bar{\omega}}_{\omega}  T_{\omega}}, \quad
{S}_{\Bar{\omega}} =\sum_{\omega \in \Omega_{SP}}{\alpha^{\Bar{\omega}}_{\omega}  S_{\omega}}.
\end{align*}
\end{definition}
We can guarantee the same objective value of the problem \eqref{model:MIP} and improve the lower bound on $\theta_\omega$ for $\omega \in \Omega_{SP}$, as shown in \cite{PartialBenders}, by including artificial scenarios and defining the second-stage decision variables for an artificial scenario as $f_{\omega}^\intercal  z_{\Bar{\omega}} = \sum_{\omega \in\Omega_{SP}}\alpha^{\Bar{\omega}}_{\omega}\theta_\omega$ in the MP.} Convex combinations of scenarios, as suggested by \cite{PartialBenders}, can create dominance relationships among scenarios, resulting in fewer feasibility cuts and a considerable reduction in the optimality gap upon termination. This technique also improves the number of instances that can be solved optimally within the time limit. Overall, including a set of scenarios $\Omega_{MP}$ in the master problem can strengthen it and lead to faster convergence. { They disclose several methods to create the artificially generated scenario $\bar \omega $, for this application, we choose to continue with choosing random values of  $\alpha^{\Bar{\omega}}_{\omega}$. }

\subsection{The TBDS Algorithm} \label{subsec:Algs}

In this section, we provide the algorithmic explanation of our proposed method. We embed the two key novelties, as explained in Sections 3.1 and 3.2, in a branch-and-cut algorithm. In addition, we use strengthened feasibility and optimality cuts as proposed by \cite{rahmaniani2017benders}.

{\SingleSpacedXI
\begin{algorithm}[!ht]
\caption{Two-Step Benders Decomposition with Scenario Clustering (TBDS)}\label{algorithm:TBDS}
\begin{algorithmic}[1]
\STATE \textbf{Input:} $AP(\cdot)$,$\{SP(\cdot,\omega)\}_{\omega \in \Omega_{SP}} $
\STATE \textbf{Output:} $(\bar x, \bar y, \Theta)$
\STATE Initialize $\Theta = -\infty$, $\theta_\omega = -\infty \quad \forall \omega \in \Omega_{SP} $, 
\State Start a callback procedure
\State \underline{Callback:} 
\IF {$(\bar x,y^*) \notin \mathcal{X}$}
    \IF {$AP(\bar x)$ is a finite number}
        \State Generate generalized Benders optimality cut \eqref{optimalitycut:generalized_benders_cut}, add to the MP
        \State Obtain $\bar y $ from $AP(\bar x)$ 
        \State Solve $SP(x^*,\Bar{y},\omega)$ and generate strengthened multi-optimality cut \eqref{optimalitycut:strmulticut}, add to the MP
    \ELSE
        \State Generate the feasibility cut \eqref{feasibilitycut}
    \ENDIF
\ELSE
     \IF {$x^*$ contains subtours}
        \State Add subtour elimination constraint of DFJ formulation
    \ENDIF
     \IF {$AP(\bar x)$ is a finite number}
         \IF {$\Theta^* < AP(\bar x)$} 
        \State Generate generalized Benders optimality cut \eqref{optimalitycut:generalized_benders_cut}, add to the MP
        \FOR {$\omega \in \Omega_{SP}$} 
        \State Obtain $\Bar{y}$ from $AP(\bar x)$
         \IF{$\theta^*_\omega < SP(\bar x,\Bar{y},\omega)$}
            \State Generate strengthened multi-optimality cut \eqref{optimalitycut:strmulticut}, add to the MP
        \ENDIF
        \ENDFOR
        \ELSE
                \State The solution $(\bar x, \Bar{y})$ is optimal
            \ENDIF
    \ELSE
        \State Generate the feasibility cut \eqref{feasibilitycut}
    \ENDIF
\ENDIF
\end{algorithmic}
\end{algorithm}}

We provide an efficient algorithmic implementation of TBDS along the lines of branch-and-cut to obtain an efficient algorithm for solving two-stage stochastic mixed-integer programs with continuous recourse. 
The TBDS algorithm dynamically adds optimality and feasibility cuts during the branch-and-bound procedure. Note the branch-and-bound procedure ensures integrality of the $x$ variables so that $x \in \mathcal X$. 

An algorithmic description of TBDS is provided in Algorithm \ref{algorithm:TBDS}. { We start with a callback function that allows us to track or modify the state of the optimization problem.} We distinct two cases for any arbitrary branch-and-bound node. First, if the associated solution $x$ is fractional (line 6), we create a feasibility cut \eqref{feasibilitycut} in case AP is infeasible (line 12). Suppose the associated solution is feasible (lines 8, 9, 10). In that case, however, we generate the strengthened optimality cuts by following Benders dual decomposition approach as these cuts are tighter than the classical optimality cuts for fractional first-stage solutions \citep{BDD}. Second, if the associated solution $x$ is integer and feasible, we { generate generalized Benders optimality cuts $\eqref{optimalitycut:generalized_benders_cut}$ and strengthened optimality cuts $\eqref{optimalitycut:strmulticut}$. Note that following the concept of Benders dual decomposition to generate strengthened optimality cuts for integer solutions is not useful as they are not as tight as the cuts following the guidelines presented in this section \citep{BDD}. { If $\bar x$ is not feasible for AP, we generate and add feasibility cut \eqref{feasibilitycut}.} 

Finally, we would like to stress that any other row generation procedure on feasibility on $\mathcal X$ and $\mathcal Y$ can easily be included in this procedure. For example, we dynamically generate subtour elimination constraints when solving the TWATSP-ST (line 15). }

\section{Implementation details of TBDS for TWATSP-ST} \label{Sec:Problem_Statement}

This section details how our TBDS approach can solve the TWATSP-ST. As stated in Section \ref{Sec:TBDS}, our TBDS method consists of a Master Problem (MP) and subproblems SP and AP. We first define the master problem of TBDS and the associated cuts in the remainder of this section. 

\subsection*{Master Program}
We state the MP, including all the associated cuts we are going to generate dynamically, for the TWATSP-ST as:
\begin{align}
(MP) :\  \min  \quad & \sum_{(i,j)  \in A}{{d_{ij}x_{ij}}} +\sum_{\omega \in \Omega_{MP}^1}p_{\omega}\left({\sum_{j \in V^+}\phi\left( e_{j\omega} +l_{j\omega}\right) + \psi o_\omega}\right) + \Theta  \label{TWATSP_Master:obj_function} \span \\
\text{s.t.} \quad & \text{First-Stage Constraints \eqref{TWATSP:cons_flow} - \eqref{TWATSP:domain_tw}}\\
&\text{Second-Stage Constraints \eqref{TWATSP:cons_departure} - \eqref{TWATSP:domain_o} for }  (\Omega_{MP}), \\
& \Theta\geq  \sum_{j \in V^+} {\sigma \left(y_{i}^{e}-y_{i}^s\right)} +  \sum_{\omega \in \Omega_{SP}} {p_\omega \theta_\omega}, \label{TWATSP_Master:cons_theta_lower_bound} \\ 
&\text{Generalized  Benders Optimality Cuts}, &  \\
& \text{Strengthened Multi-Optimality Cuts }  (\Omega_{SP}), &  \\
&\text{Feasibility Cuts }  (\Omega_{SP}) , &  \\
& \Theta \in \mathbb{R}, \label{TWATSP_Master:cons_theta_domain} \\
& \theta_\omega \in \mathbb{R} & \forall \ \omega \in \Omega_{SP}.  \label{TWATSP_Master:cons_theta_omega_domain} 
\end{align}

{The first term of the objective function of the MP is the total distance traveled, the second term is the expected time window exceedances, and the overtime cost of scenarios included in the MP (set~$\Omega^1_{MP}$), and the last term provides a lower bound on the expected second-stage cost and the time window assignment cost. We include first-stage constraints \eqref{TWATSP:cons_flow} - \eqref{TWATSP:domain_tw} and the second-stage constraints associated with the scenario set $\Omega_{MP}$ in the MP. Constraint \eqref{TWATSP_Master:cons_theta_lower_bound} is the auxiliary variables $\Theta$ (and $\theta(\omega)$) for recourse cost of SP and AP. The optimality and feasibility cuts are made specific in the next subsection. Finally, Constraints \eqref{TWATSP_Master:cons_theta_domain} and \eqref{TWATSP_Master:cons_theta_omega_domain} restrict the domain of the auxiliary variables $\Theta$ and $\theta_\omega$.}

\subsection{Cuts}  \label{subsec:TWATSP-ST_Cuts}
In line with Section \ref{Sec:TBDS}, we define {an aggregated subproblem AP} taking as input the binary first-stage decisions, i.e., the routing decision in the TWATSP-ST. The second set of subproblems SP for each scenario $\omega \in \Omega_{SP}$ are then obtained via the two-step decomposition as outlined in Section \ref{Sec:TBDS}, i.e., as input we take the binary first-stage decisions and the solution of the time window assignment variables after solving AP. We state for the TWATSP-ST the generalized Benders optimality cut \eqref{optimalitycut:generalized_benders_cut} introduced in Section \ref{subsec:TwoStepDecomposition} as follows. 
\subsubsection*{Generalized Benders Optimality Cut.}
For any solution $\hat x$ of MP-2 (or MP) that is feasible for AP, let $\bar \lambda$ indicate the values of the dual multipliers of the constraint $\{x = \hat x\}$ in $AP(\hat x)$ and  (${\bar x, \bar y, \bar z}$) be the optimal solution of $AP(\hat x)$. The Generalized Benders Optimality cut is then given by: 
\begin{equation} \label{optimalityCut:TWATSP_generalizedBendersCut}
\begin{split}
{\Theta \geq \sum_{j \in V^+} {\sigma \left(\bar y_{i}^e- \bar y_{i}^s\right)} + \sum_{\omega \in \Omega_{SP}} {p_\omega\left(\sum_{j \in V^+}\phi\left( \bar e_{j\omega} +  \bar l_{j\omega}\right) +\psi \bar o_{\omega}  \right) + \sum_{(i,j) \in A}\bar \lambda_{ij}(x_{ij}-\bar x_{ij}) }}
\end{split}
\end{equation}
For the second step of our decomposition, we define $SP(\bar x, \bar y, \omega)$ and create strengthened multi-optimality cuts \eqref{optimalitycut:strmulticut} as detailed in Section \ref{subsec:TwoStepDecomposition}.
\subsubsection*{Strengthened Multi-Optimality Cuts.}
For any solution $\hat x$ of MP-2 (or MP) that is feasible for AP and the solution $(\bar y = \{\bar y^{s},\ \bar y^{e}\})$ of $AP(\hat x)$ for $\omega \in \Omega_{SP}$, let ($\tilde x, \tilde y, \tilde z$) be the optimal solution of $SP(\hat x, \bar y, \omega)$ and  ($\tilde \nu $, $\tilde \eta^{s} $, $\tilde \eta^{e}$) indicate the values of the dual multipliers related to constraints $\{x = \hat x$, $y^{s} = \bar y^{s}$, $y^{e} = \bar y^{e} \}$  respectively; then, 
\begin{equation}
{\theta_\omega \geq \sum_{j \in V^+}\left(\phi\left(\tilde{e}_{j\omega} +\tilde{l}_{j\omega} + \psi \tilde{o}_{\omega}\right)+ \tilde \eta^{s}_j(y_j^{s}- \tilde y_j^{s}) + \tilde \eta^{e}_j(y_j^{e}- \tilde y_j^{e}) \right) +  \sum_{(i,j) \in A}\tilde \nu_{ij}(x_{ij}-\tilde x_{ij}) } \label{optimalitycut:TWATSP_multicut}   
\end{equation} 
is the strengthened multi-optimality cut for $\omega \in \Omega_{SP}$ for MP. We emphasize that SP takes the routing decisions from MP and time window assignments from AP as input. Via the auxiliary decision variable $\theta_\omega$, we approximate the cost function of each scenario $\omega \in  \Omega_{SP}$ (or the objective function value of each SP) with lifted strengthened multi-optimality cuts. 

As a final remark, we note that TWATSP-ST has the relatively complete recourse property. In other words, the second-stage decision is always feasible for any feasible first-stage decision (routing and time window variables). Therefore, we do not need to check feasibility in the branch-and-cut tree when $\bar x$ is integer. However, for a fractional node, we state the following feasibility cut.
\subsubsection*{Feasibility Cuts.} 
If the solution $\hat x$ of MP-2 (or MP) is infeasible for AP, we generate the feasibility cut
\begin{equation*}
    0 \geq \sum_{i \in V^+} \breve\epsilon + \sum_{(i,j) \in A} \breve\lambda_{ij}(x_{ij}-\breve x_{ij}) \label{feasibilitycut:TWATSP}
\end{equation*}
where $\breve\epsilon$  and $\breve x $ are the optimal values of the $\epsilon$ and $z$ variables in the feasibility problem and we set $\breve \lambda$ as the value of the dual variable associated with the constraint $\{ x=\hat x \}$.

\section{Computational Results} \label{Sec:Comp_Results}
{ In this section, we present the computational results of solving the TWATSP-ST with our new TBDS method. We first determine the impact of our proposed strategies to strengthen the cut generation process of Benders decomposition. We then seek to understand the source of the computational benefits associated with using TBDS. { Finally, we present an analysis of the problem context.}} We develop six variants of TBDS to assess and structurally benchmark the two main novelties of TBDS: the two-step decomposition to generate { \textit{strengthened} multi and generalized Benders optimality cuts} and the new scenario-retention strategy to include scenarios in master problem. The six considered variants are:

\begin{enumerate}
\item \textbf{BD} uses standard \textit{multi-optimality cuts} utilizing Benders dual decomposition, as introduced by \cite{BDD}.
\item \textbf{TBD} extends BD by including the decomposition over the integer and linear first-stage variables. This variant tests the impact of our { strengthened multi-optimality cuts and generalized Benders optimality cuts} compared with Benders dual decomposition. 
\item \textbf{BDP} extends BD by including the first-stage constraints on the master problem of randomly chosen and artificial scenarios. This variant combines partial Benders decomposition \citep{PartialBenders} with Benders dual decomposition.
\item \textbf{TBDP} extends BDP by { strengthened multi-optimality cuts}. In other words, this is our TBDS method without the new scenario-retention strategy.
\item \textbf{BDS} extends BDP by the new scenario-retention strategy. Thus, this method combines benders dual decomposition and partial benders decomposition with our new scenario retention strategy but does not include our two-step decomposition idea and the associated strengthened multi-optimality cuts.
\item \textbf{TBDS} extends BDS by including the two-step first-stage decomposition. This is our TBDS method as presented in Sections \ref{Sec:TBDS} and \ref{Sec:Problem_Statement}, and thus includes the two-step decomposition and the scenario-retention strategy.
\end{enumerate}
{The remainder of this section is structured as follows.
We introduce a new set of benchmark instances in Section \ref{subsec:TestInstances}. Details on the parameter tuning are provided in Section \ref{subsec:ParameterSetting}. 
We compare the performance of the six variants of TBDS on the new benchmark instances in Section \ref{subsec:ResultsBC}. Further details explaining the effectiveness of TBDS are given in Section \ref{subsec:detailsofperformance}.}

\subsection{Benchmark Instances} \label{subsec:TestInstances}

We adapt benchmark instances from the literature to our problem setting, generating 126 new benchmark instances. Specifically, we derive 56 instances with clustered customer locations based on single vehicle routes from Solomon's VRPTW-RC instances proposed by \cite{potvin1996vehicle} and 70 instances from \cite{gendreau1998generalized}, which include not only the customer locations but also the associated service times. We refer to these as \textit{rc\_} and \textit{n\_w\_} instances, respectively. The end of the depot's time window defines the shift length $T$. All instances and associated solutions are included as supplementary material to {this paper}.

To balance the four penalty terms in the objective function and make them directly comparable concerning the routing cost for most problem instances, we used penalty weights of $\phi=3$ and $\sigma=1$ to penalize expected delay/earliness and time window width at customer $i \in V^+$, respectively. Comparably to \cite{jabali2015self}, we set the penalty for expected shift overtime to $\psi=4$.

We conduct a Sample Average Approximation analysis to determine the number of scenarios to correctly represent uncertainty, considering our routing solutions' stability as the number of customers in the instances increases. This results in 100 scenarios, which we use throughout all experiments. 

We follow a similar approach for generating random travel times as \cite{jabali2015self} and \cite{vareias2019assessing}. For each scenario, we determine the travel time $t_{ij}$ for arc $(i,j) \in A$  by adding a random disruption parameter $\delta_{ij}$ to the Euclidean distance $d_{ij}$. We assume a Gamma-distributed disruption parameter with shape $k$ and scale $\theta_{ij}$ depending on the distance and a coefficient of variation $cov = 0.25$ to model realistic travel times and disruptions. We define $\eta = 0.35$ as a congestion level, representing the expected increase in travel time after a disruption occurs, i.e., $\mathbb{E}[\delta_{ij}] = \eta d_{ij}$. We assume $\delta_{ij} \sim G(k,\theta_{ij}) $, resulting in
\begin{align*}
    \mathbb{E}[\delta_{ij}] &= k\theta_{ij} = \eta d_{ij}, \\
   Var(\delta_{ij}) & = k\theta_{ij}^2. 
\end{align*}
Here, parameters $k$ and $\theta_{ij}$ are set according to
\begin{equation*}
    \quad k = \dfrac{1}{{cov}^2}, \quad\theta_{ij} = \eta d_{ij}{cov}^2
\end{equation*}

\subsection{Parameter Settings and Implementation Details} \label{subsec:ParameterSetting}

The master problem \eqref{TWATSP_Master:obj_function}-\eqref{TWATSP_Master:cons_theta_omega_domain} contains both a set of selected scenarios $\Omega_{MP}^1$ and some artificial scenarios $\Omega_{MP}^2$. Section \ref{subsec:ScenarioSelection} explains in detail how to select and create artificial scenarios. The performance of our method depends on the number of scenarios added to the master problem. {By varying the fraction of scenarios added to the master problem for actual scenarios (5\%, 10\%, 15\%) and artificial scenarios (5\%, 10\%),} we evaluate the performance of the methods BDP, TBDP, BDS, and TBDS under different scenario selection strategies. {To create artificial scenarios, we randomly select the values for the convex composition}. The results presented for the methods above are obtained by selecting the best combination of scenario levels corresponding to the highest number of solved instances in each method. 

We embed all method variants in a branch-and-cut framework. We only add Benders cuts \eqref{optimalityCut:TWATSP_generalizedBendersCut}-\eqref{optimalitycut:TWATSP_multicut} at the root node and upon discovery of a possible incumbent solution.
Subtour elimination constraints are included dynamically. All algorithms are coded in Python 3.9.7 in combination with Gurobi 9.5.1. The experiments were run on a virtual machine with 32 CPU cores and 64 GB RAM under the Linux operating system, which was sufficient for all experiments. All algorithms have a time limit of three hours. All computation times reported are in minutes. 

\subsection{Overall performance of TBDS}
\label{subsec:ResultsBC}
\begin{figure}[!htbp]
\centering
\caption{Percentage of Instances Solved Within 3 hours}
\includegraphics[width=0.6\linewidth]{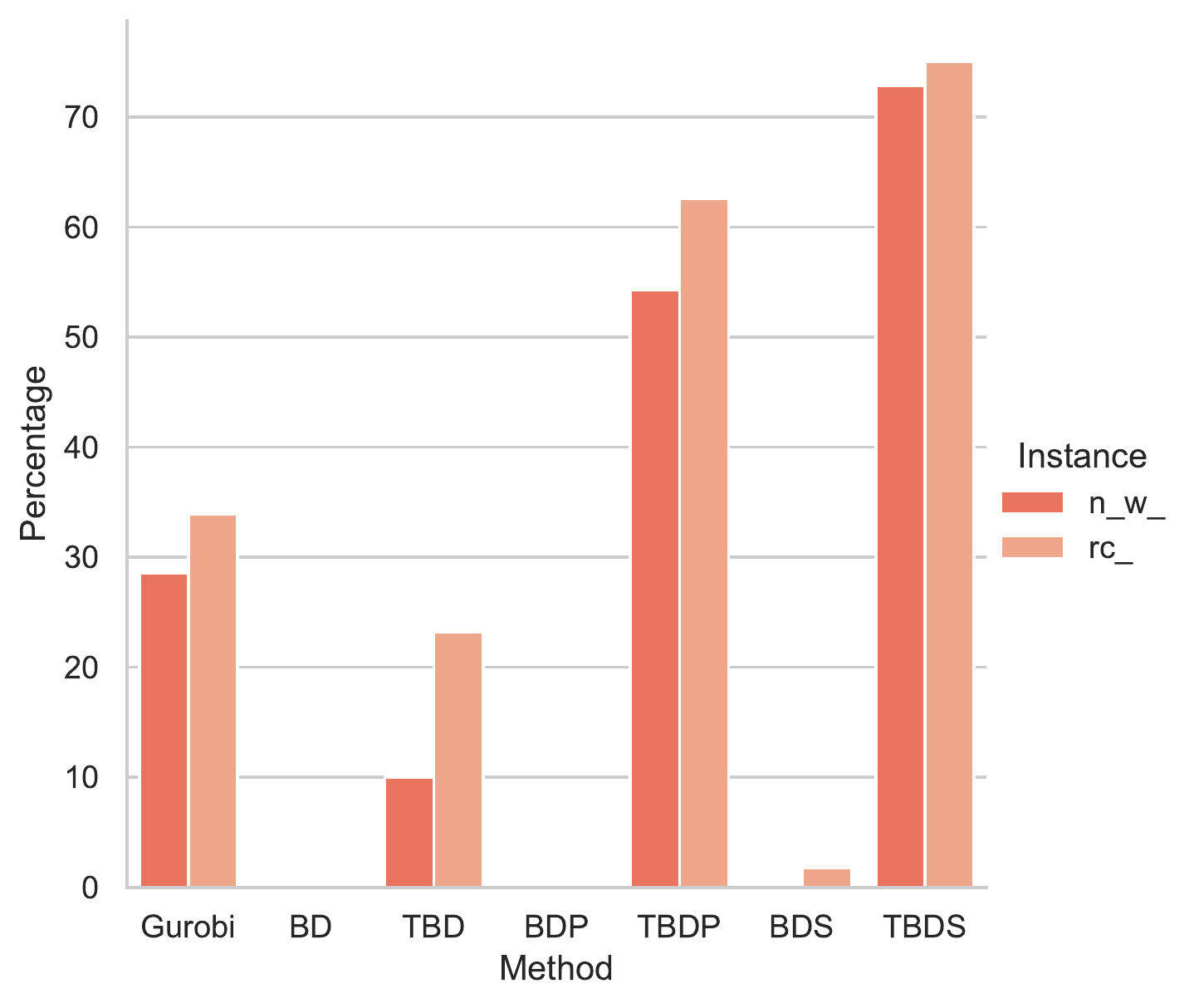}
\label{plot:solvedinstances}
\end{figure}

We start by analyzing the ability of the six method variants and of a straightforward Gurobi implementation to solve our benchmark instances. Figure \ref{plot:solvedinstances} illustrates the fraction of solved instances within a three-hour time limit. In Tables \ref{tab:Performance_nw} and \ref{tab:Performance_rc}, we present the average optimality gap at the termination for each method variant and the average amount of time until termination. 

These results show that our proposed solution approach, {TBDS}, solves up to 25 customers in both sets of benchmark instances. Looking at the total number of solved instances within the time limit displayed in Figure \ref{plot:solvedinstances}, it is evident that the variants with two-step decomposition outperform the variants without two-step decomposition. Specifically, TBDS solved 72.8\% and 75.2\% of instances in both benchmark sets to optimality, while BDS only solved 1.78\% of instances. In contrast, BD and BDP failed to solve any instance. The TBD and TBDP variants were able to solve between 10\% and 62.5\% of all benchmark instances. The significant increase in the performance of TBDS can be attributed to the use of a new scenario-retention strategy and the two-step decomposition over continuous and binary first-stage variables.

\begin{table}[!htbp]
\small
\centering
\caption{Performance of the six TBDS variants on n\_w\_ instances}
\label{tab:Performance_nw}
\begin{tabular}{@{}llrrrrrr@{}}
\toprule
Cust.           &                                     & BD    & BDP    & BDS    & TBD    & TBDP   & TBDS   \\ \midrule
\multirow{3}{*}{10} & {\#   Solved (/10)}   & 0     & 0     & 0     & 4     & 10    & 10    \\
                    & {Time to opt. (min.)} & 180.0 & 180.0 & 180.0 & 45.6  & 4.2   & 13.1  \\
                    & Optimality gap                    & 94.5  & 93.4  & 93.5  & 20.2  & 0.0   & 0.0   \\ \midrule
\multirow{3}{*}{13} & {\# Solved (/10)}     & 0     & 0     & 0     & 2     & 8     & 9     \\
                    & {Time to opt. (min.)} & 180.0 & 180.0 & 180.0 & 67.9  & 6.0   & 12.7  \\
                    & Optimality gap                    & 93.9  & 91.4  & 91.7  & 20.4  & 24.1  & 16.2  \\\midrule
\multirow{3}{*}{15} & {\# Solved (/10)}     & 0     & 0     & 0     & 1     & 10    & 10    \\
                    & {Time to opt. (min.)} & 180.0 & 180.0 & 180.0 & 75.6  & 18.0  & 15.7  \\
                    & Optimality gap                    & 93.6  & 91.0  & 91.0  & 23.9  & 0.0   & 0.0   \\\midrule
\multirow{3}{*}{18} & {\# Solved (/10)}     & 0     & 0     & 0     & 0     & 8    & 10    \\
                    & {Time to opt. (min.)} & 180.0 & 180.0 & 180.0 & 180.0 & 28.4  & 18.5  \\
                    & Optimality gap                    & 92.7  & 88.8  & 89.0  & 26.1  & 0.0   & 0.0   \\\midrule
\multirow{3}{*}{20} & {\# Solved (/10)}     & 0     & 0     & 0     & 0     & 2    & 10    \\
                    & {Time to opt. (min.)} & 180.0 & 180.0 & 180.0 & 180.0 & 20.4  & 23.1  \\
                    & Optimality gap                    & 92.3  & 89.3  & 89.2  & 28.3  & 0.0   & 0.0   \\\midrule
\multirow{3}{*}{23} & {\# Solved (/10)}     & 0     & 0     & 0     & 0     & 0     & 1     \\
                    & {Time to opt. (min.)} & 180.0 & 180.0 & 180.0 & 180.0 & 180.0 & 18.5  \\
                    & Optimality gap                    & 91.7  & 85.8  & 86.8  & 34.3  & 31.9  & 33.6  \\\midrule
\multirow{3}{*}{25} & {\# Solved (/10)}     & 0     & 0     & 0     & 0     & 0     & 1     \\
                    & {Time to opt. (min.)} & 180.0 & 180.0 & 180.0 & 180.0 & 180.0 & 135.3 \\
                    & Optimality gap                    & 91.2  & 86.6  & 84.7  & 47.6  & 45.9  & 39.4 \\ \bottomrule
\end{tabular}
\end{table}
\begin{table}[!htbp]
\small
\centering
\caption{Performance of the six TBDS variants on rc\_ instances}
\label{tab:Performance_rc}
\begin{tabular}{@{}llrrrrrr@{}}
\toprule
Cust. &                                                 & BD    & BDP    & BDS    & TBD    & TBDP   & TBDS \\ \midrule
\multirow{3}{*}{10} & {\#   Solved (/10)}   & 0     & 0     & 1     & 3     & 8     & 8     \\
                    & {Time to opt. (min.)} & 180.0 & 180.0 & 180.0 & 21.0  & 3.1   & 10.7  \\
                    & Optimality gap                    & 92.0  & 90.8  & 90.6  & 11.6  & 0.0   & 0.0   \\\midrule
\multirow{3}{*}{13} & {\# Solved (/10)}     & 0     & 0     & 0     & 3     & 8     & 8     \\
                    & {Time to opt. (min.)} & 180.0 & 180.0 & 180.0 & 58.7  & 8.7   & 20.6  \\
                    & Optimality gap                    & 91.3  & 88.8  & 89.0  & 14.4  & 0.0   & 0.0   \\\midrule
\multirow{3}{*}{15} & {\# Solved (/10)}     & 0     & 0     & 0     & 3     & 8     & 8     \\
                    & {Time to opt. (min.)} & 180.0 & 180.0 & 180.0 & 98.5  & 13.1  & 23.1  \\
                    & Optimality gap                    & 90.7  & 88.7  & 87.3  & 16.9  & 0.0   & 0.0   \\\midrule
\multirow{3}{*}{18} & {\# Solved (/10)}     & 0     & 0     & 0     & 3     & 7     & 8     \\
                    & {Time to opt. (min.)} & 180.0 & 180.0 & 180.0 & 180.0 & 25.8  & 82.7  \\
                    & Optimality gap                    & 89.7  & 86.8  & 86.5  & 19.7  & 0.0   & 0.0   \\\midrule
\multirow{3}{*}{20} & {\# Solved (/10)}     & 0     & 0     & 0     & 1     & 4     & 8     \\
                    & {Time to opt. (min.)} & 180.0 & 180.0 & 180.0 & 188.1 & 26.0  & 32.0  \\
                    & Optimality gap                    & 89.0  & 85.2  & 85.5  & 17.6  & 0.0   & 0.0   \\\midrule
\multirow{3}{*}{23} & {\# Solved (/10)}     & 0     & 0     & 0     & 0     & 0     & 0     \\
                    & {Time to opt. (min.)} & 180.0 & 180.0 & 180.0 & 180.0 & 180.0 & 180.0 \\
                    & Optimality gap                    & 87.9  & 84.8  & 84.5  & 24.0  & 25.9  & 21.3  \\\midrule
\multirow{3}{*}{25} & {\# Solved (/10)}     & 0     & 0     & 0     & 0     & 0     & 2     \\
                    & {Time to opt. (min.)} & 180.0 & 180.0 & 180.0 & 180.0 & 180.0 & 62.9  \\
                    & Optimality gap                    & 87.2  & 83.1  & 84.2  & 26.5  & 22.4  & 23.1 \\ \bottomrule
\end{tabular}
\end{table}

Next, we focus on looking at the solution time till termination and the average optimality gap if the instance is unsolved. Tables \ref{tab:Performance_nw} and \ref{tab:Performance_rc} provide detailed results for each benchmark set. Noting that the BD and BDP variants fail to solve every instance within the 3-hour time limit, the computational time of the solved instance in the BDS is close to 180 minutes. We see a clear decrease in computational time ranging between 17.7 and 65.7 minutes in the variants with the two-step decomposition {with the largest share of the runtime stemming from the master problem}. Comparing BD, BDP, and BDS, the average optimality gap of non-solved instances decreases significantly. Among variants without two-step decomposition, the average optimality gap decreases from 91.1\% to 88.2\%.

Similarly, with the variants with two-step decomposition, we see that the average optimality gap decreases from 31.5\% to 25.5\%. The results clearly show that our TBDS outperforms the state-of-the-art, as it solves more instances to optimality and obtains smaller optimality gaps for the instances unable to solve optimally. We can say that within the time limit, for the instances that benchmark methods are unable to solve,  TBDS produces a solution of provably higher quality than the benchmark methods. Essential is combining the two main ideas of TBDS, i.e., the two-step decomposition and the new scenario-retention strategy.

{
\subsection{Details on the performance of TBDS} \label{subsec:detailsofperformance}

In this section, we answer multiple questions regarding the methods' performances. We begin by analyzing the convergence rate of the benchmark algorithms. Afterward, we focus on the impact of our strengthened multi-optimality cuts and our new scenario-retention strategy by analyzing the root node optimality gaps. Recall that we solve the TWATSP-ST with the aforementioned cuts with a branch-and-cut algorithm wherein cuts are generated throughout the tree search. Since we did not observe a significant difference between the two benchmark sets, we present the results for convergence rate for both instance sets without differentiation. 

\subsubsection{Convergence Rate}
Figure \ref{fig:convergencerate} compares the number of iterations required (on average) and the average number of open and explored nodes in the tree search at termination. We see that using two-step decomposition either on its own or combined with a new scenario-retention strategy greatly speeds up convergence. We also note that the new scenario-retention strategy lowers the number of iterations needed rather than random selection or no selection of scenarios. 
We next look into the size of the branch-and-bound tree. We consider the average number of open and explored nodes in the branch-and-bound tree at termination to measure convergence. Figure \ref{fig:nodes} shows the average number of nodes at the termination for each method. Including strengthened multi and generalized Benders, optimality cuts result in much smaller search trees. These results also support the findings that TBDS yields high-quality bounds at the root node, which will be analyzed in detail in the next section. 

\begin{figure}[!htbp]
\centering
\caption{Convergence rate by Method}
\begin{subfigure}[b]{0.40\textwidth}
     \centering
     \includegraphics[width=\textwidth]{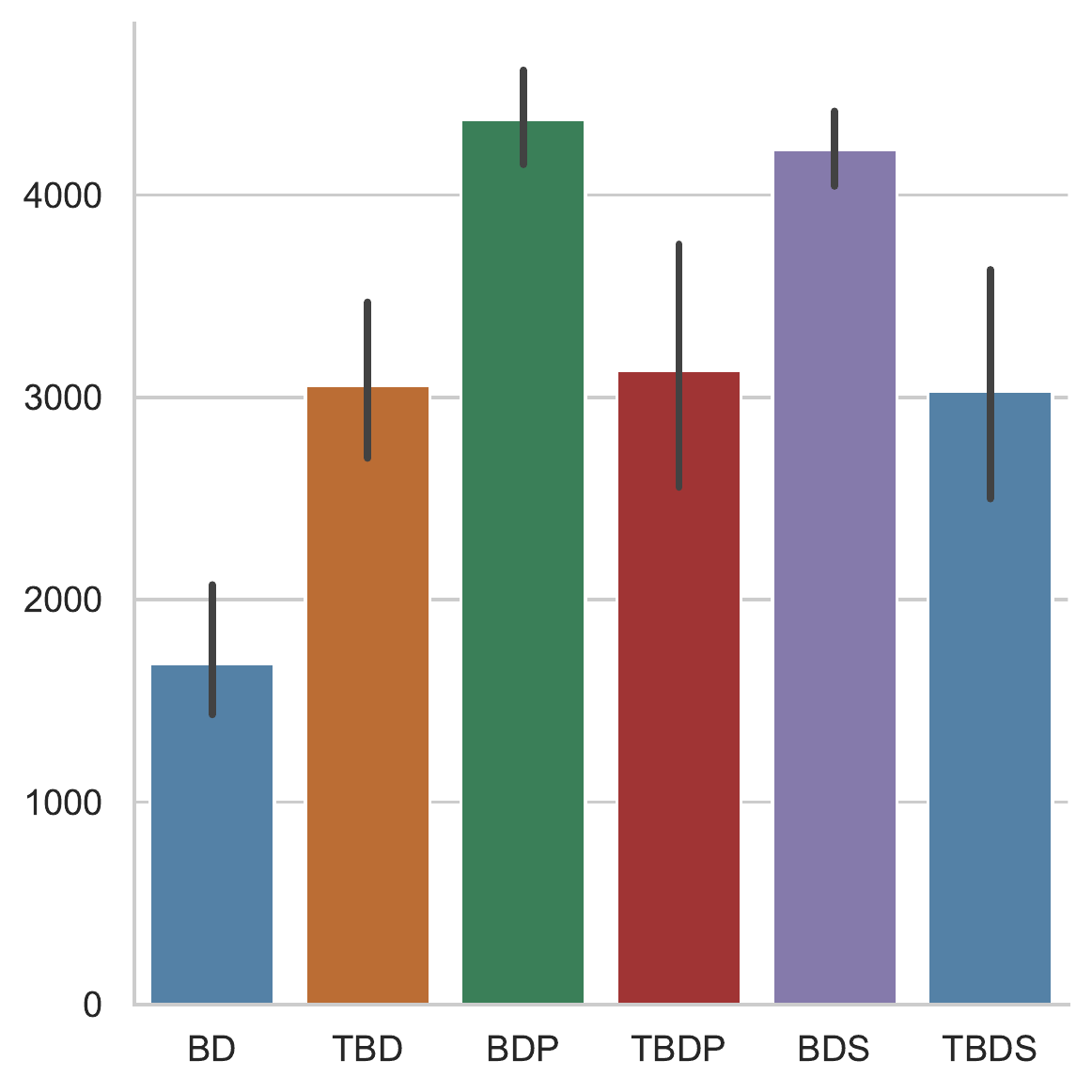}
     \caption{Number of iterations}
     \label{fig:iteration}
 \end{subfigure}
\hfill
\begin{subfigure}[b]{0.47\textwidth}
    \centering                
    \includegraphics[width=\textwidth]{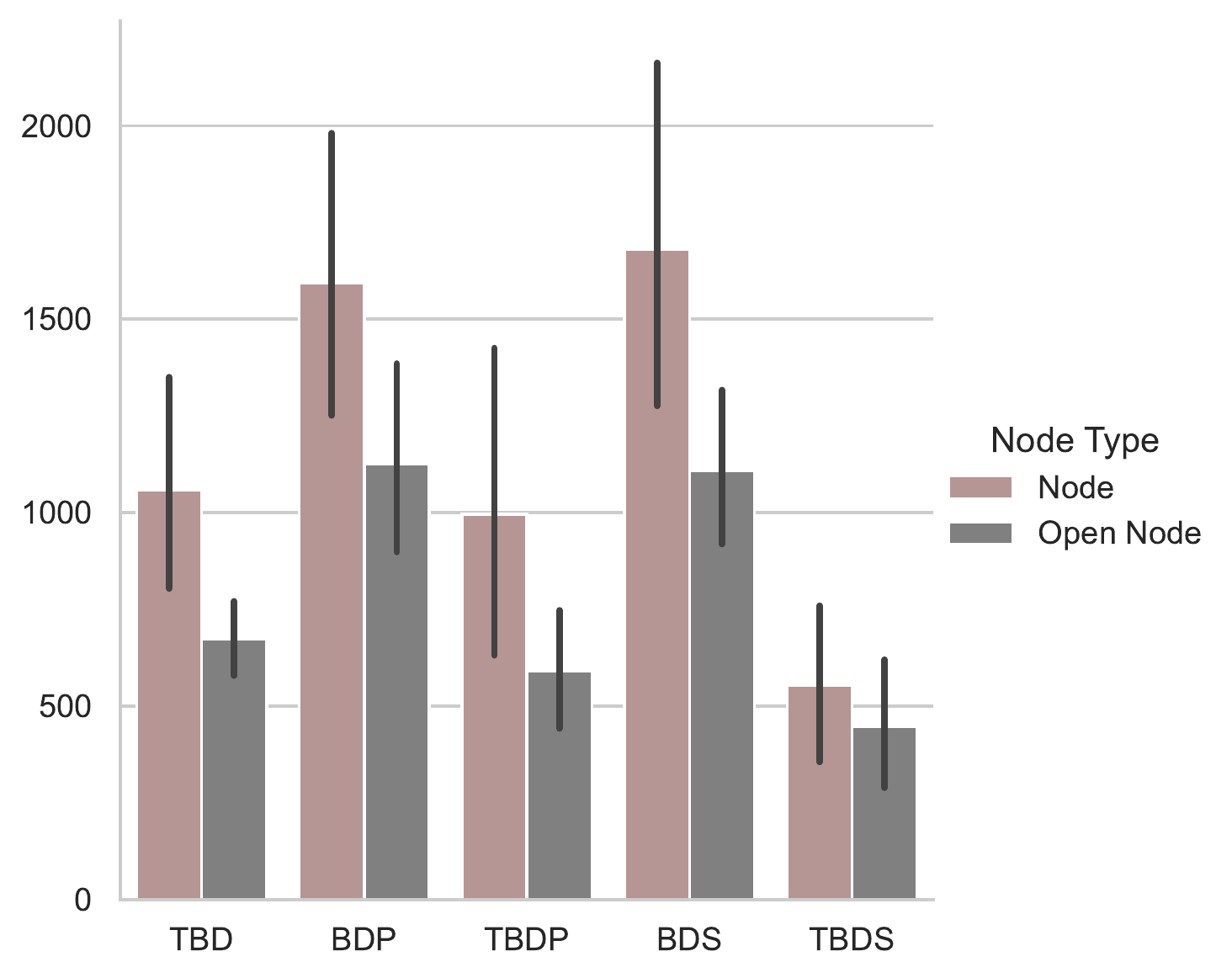}
    \caption{Number of nodes in the search tree}
    \label{fig:nodes}
\end{subfigure}
\label{fig:convergencerate}
\end{figure}

Another statistic to measure convergence is the number of Benders cuts generated. As such, Figure \ref{fig:Cuts} reports the number of optimality and feasibility cuts generated for different methods. We consider both the cuts generated when solving the linear programming relaxation of the TWATSP-ST and solving the TWATSP-ST via branch-and-cut. 
\begin{figure}[!htbp]
\centering
\caption{Number of optimality and feasibility cuts by method}
\begin{subfigure}[b]{0.3\textwidth}
     \centering
     \includegraphics[width=\textwidth]{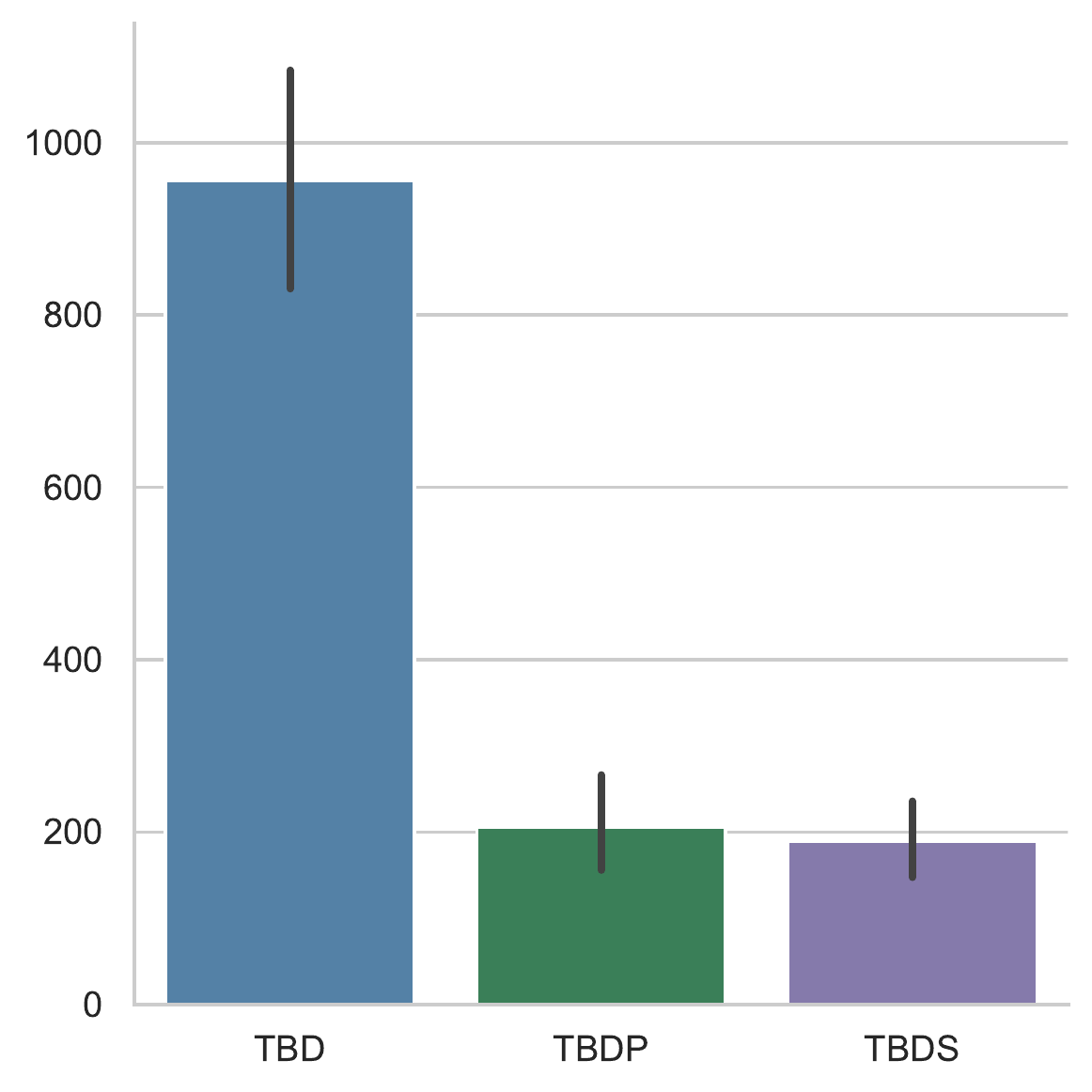}
     \caption{Generalized cuts}
     \label{fig:generalizedcut}
 \end{subfigure}
\hfill
\begin{subfigure}[b]{0.3\textwidth}
    \centering                
    \includegraphics[width=\textwidth]{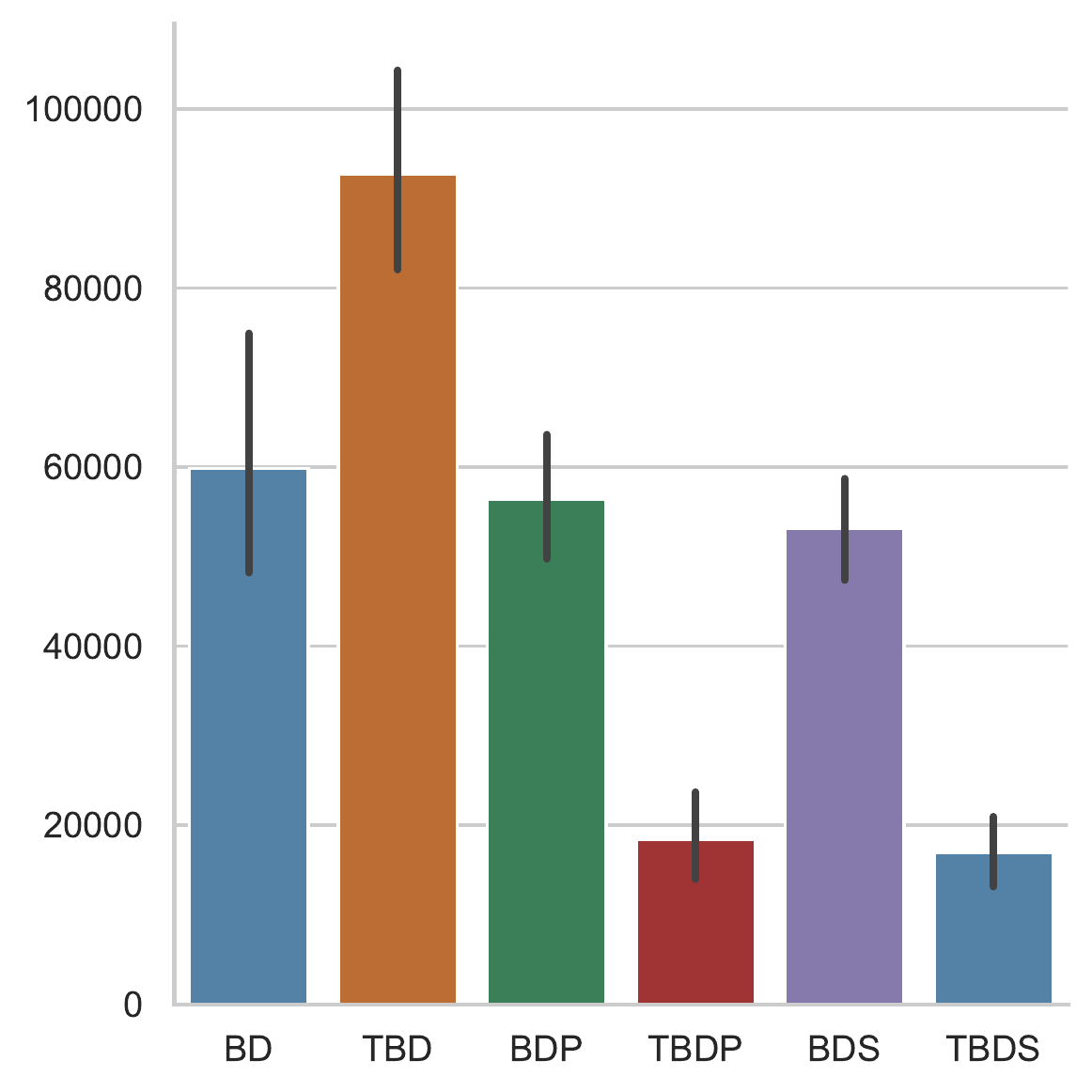}
    \caption{Multi-Optimality cuts}
    \label{fig:multioptimalitycut}
\end{subfigure}
\hfill
\begin{subfigure}[b]{0.3\textwidth}
    \centering                
    \includegraphics[width=\textwidth]{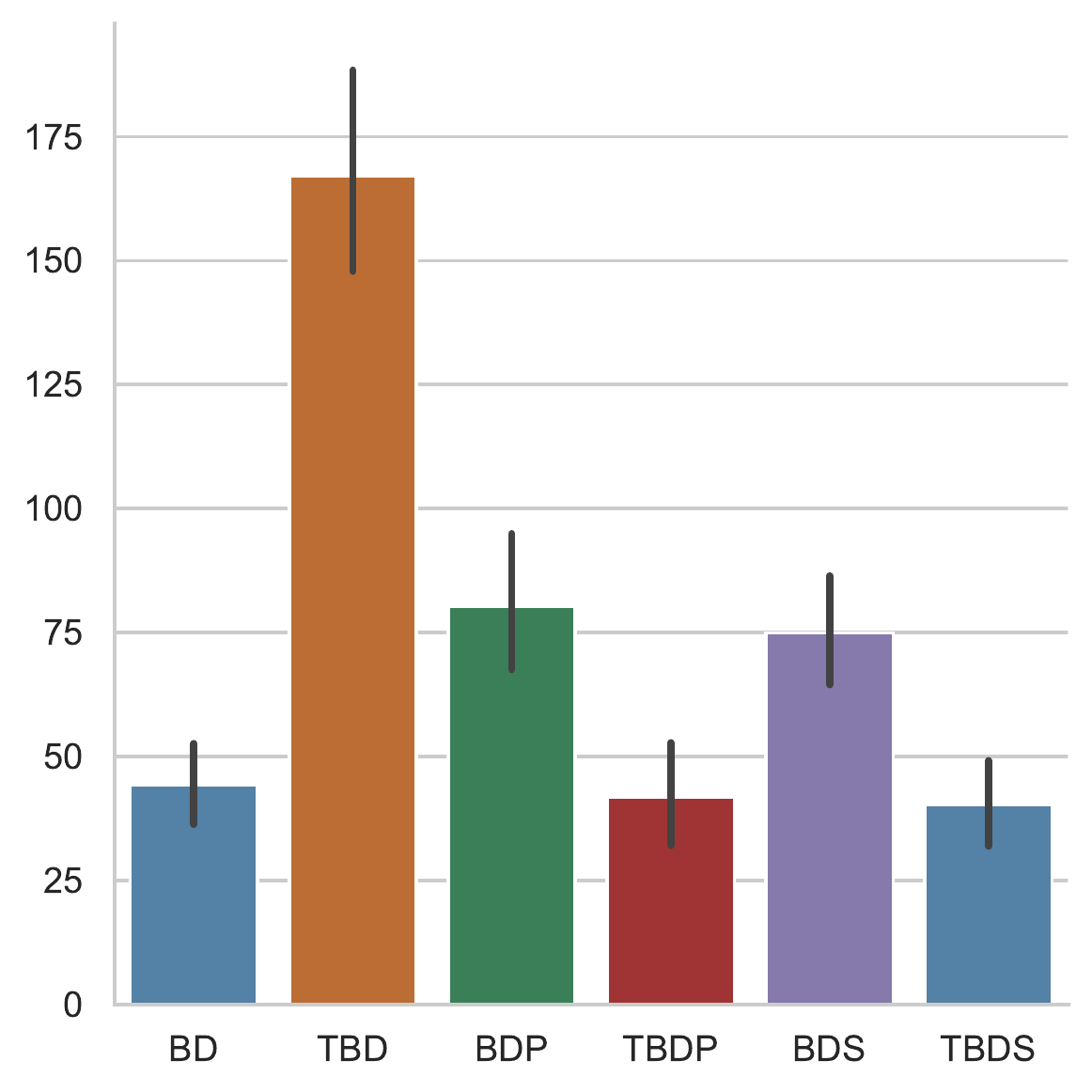}
    \caption{Feasibility cuts}
    \label{fig:feascuts}
\end{subfigure}
\label{fig:Cuts}
\end{figure}

We observe from the results that the new scenario-retention strategy and having a two-step decomposition yields a significant decrease in optimality cuts (both generalized and multi-cuts) and feasibility cuts. Two-step decomposition significantly reduces the amount of feasibility cuts and optimality cuts needed. We look into the bounds found at the root node to understand why two-step decomposition and new scenario-retention strategy impact the cuts generated throughout the search. We discuss these results in the following section. 

}
\subsubsection{Bounds at the Root Node}
\label{subsec:Rootnode}
Table \ref{tab:Bounds_nw} and \ref{tab:Bounds_rc} compare root-node lower bounds (\textit{LB}), root-node upper bounds (\textit{UB}), and root-node gaps of the six TBDS variants (\textit{BD, TBD, BDP, TBDP, BDS, TBDS}) for different instance sizes ({10, 13, 15, 18, 20, 23, 25}), listed in the first column. 
The \textit{Root Node Gap} is computed as $((UB-LB)/LB)\times 100)$ and gives the relative difference between the lower and upper bound at the root node. The lower and upper bound information is obtained as the root node in evaluating the full branch-and-cut algorithm after the first branching decision, which partly depends on the procedures embedded within Gurobi.

\begin{table}[!htbp]
\caption{Average Root Node Lower and Upper Bounds of n\_w\_ Instances for the TBDS Variants}
\label{tab:Bounds_nw}
\small
\centering
\begin{tabular}{@{}llrrrrrr@{}}
\toprule
\# Customers & & BD & BDP & BDS & TBD & TBDP & TBDS \\ \midrule
\multirow{3}{*}{10}    & Lower Bound        & 407.9   & 332.7  & 339.5  & 248.2  & 322.6  & 331.6  \\
                       & Upper Bound        & 5586.1  & 4778.7 & 4775.0 & 314.9  & 386.7  & 372.0  \\
                       & Root Node Gap (\%) & 92.7    & 93.0   & 92.9   & 21.2   & 14.2   & 13.3   \\ \midrule
\multirow{3}{*}{13}    & Lower Bound        & 516.5   & 446.7  & 439.0  & 518.6  & 446.7  & 441.5  \\
                       & Upper Bound        & 9914.3  & 4892.7 & 4874.5 & 634.0  & 510.3  & 512.9  \\
                       & Root Node Gap (\%) & 94.8    & 90.9   & 91.0   & 18.2   & 12.5   & 13.9  \\ \midrule 
\multirow{3}{*}{15}    & Lower Bound        & 571.6   & 536.0  & 555.3  & 335.3  & 537.8  & 519.2  \\
                       & Upper Bound        & 13446.9 & 4982.0 & 4990.8 & 440.0  & 634.1  & 608.4  \\
                       & Root Node Gap (\%) & 95.7    & 89.2   & 88.9   & 23.8   & 15.2   & 14.7   \\ \midrule
\multirow{3}{*}{18}    & Lower Bound        & 689.5   & 590.6  & 627.1  & 323.5  & 577.4  & 635.4  \\
                       & Upper Bound        & 18862.6 & 5036.6 & 5062.6 & 428.7  & 692.3  & 775.3  \\
                       & Root Node Gap (\%) & 96.3    & 88.3   & 87.6   & 24.5   & 16.6   & 18.0   \\ \midrule
\multirow{3}{*}{20}    & Lower Bound        & 753.7   & 659.0  & 678.5  & 353.9  & 667.4  & 663.4  \\
                       & Upper Bound        & 23932.4 & 5105.0 & 5114.0 & 480.0  & 807.6  & 804.8  \\
                       & Root Node Gap (\%) & 96.9    & 87.1   & 86.7   & 26.3   & 17.4   & 17.6   \\ \midrule
\multirow{3}{*}{23}    & Lower Bound        & 845.4   & 831.1  & 800.0  & 845.4  & 831.1  & 800.0  \\
                       & Upper Bound        & 41065.3 & 5277.1 & 5235.5 & 1086.2 & 1031.1 & 994.2  \\
                       & Root Node Gap (\%) & 97.9    & 84.3   & 84.7   & 22.2   & 19.4   & 19.5   \\ \midrule
\multirow{3}{*}{25}    & Lower Bound        & 923.5   & 957.5  & 902.0  & 380.6  & 961.7  & 928.4  \\
                       & Upper Bound        & 43749.3 & 5403.5 & 5337.5 & 535.1  & 1610.9 & 1591.8 \\
                       & Root Node Gap (\%) & 97.9    & 82.3   & 83.1   & 28.9   & 40.3   & 41.7   \\\bottomrule
\end{tabular}
\end{table}
\begin{table}[!htbp]
\caption{Average Root Node Lower and Upper Bounds of rc\_ Instances for All TBDS Variants}
\label{tab:Bounds_rc}
\small
\centering
\begin{tabular}{@{}llrrrrrr@{}}
\toprule
\# Customers &                               & BD    & BDP      & BDS     & TBD     & TBDP    & TBDS \\ \midrule
\multirow{3}{*}{10}    & Lower   Bound      & 676.4   & 597.9  & 648.8  & 397.2  & 626.9  & 647.2  \\
                       & Upper Bound        & 5988.8  & 5043.9 & 5147.0 & 457.8  & 718.2  & 690.6  \\
                       & Root Node Gap (\%) & 88.7    & 88.1   & 87.4   & 13.2   & 9.8    & 9.3    \\ \midrule
\multirow{3}{*}{13}    & Lower Bound        & 756.6   & 729.6  & 730.1  & 756.6  & 683.9  & 751.7  \\
                       & Upper Bound        & 5706.6  & 5175.6 & 5228.3 & 870.2  & 778.2  & 831.6  \\
                       & Root Node Gap (\%) & 86.7    & 85.9   & 86.0   & 13.1   & 12.1   & 9.6    \\ \midrule
\multirow{3}{*}{15}    & Lower Bound        & 769.9   & 810.3  & 809.5  & 421.7  & 802.1  & 815.0  \\
                       & Upper Bound        & 12160.6 & 5256.3 & 5307.6 & 473.9  & 929.8  & 923.1  \\
                       & Root Node Gap (\%) & 93.7    & 84.6   & 84.7   & 11.0   & 12.3   & 13.1   \\ \midrule
\multirow{3}{*}{18}    & Lower Bound        & 999.9   & 996.3  & 1004.0 & 597.4  & 940.1  & 1002.2 \\
                       & Upper Bound        & 19988.2 & 5442.3 & 5502.1 & 718.4  & 1097.1 & 1168.8 \\
                       & Root Node Gap (\%) & 95.0    & 81.7   & 81.8   & 16.9   & 14.3   & 14.3   \\ \midrule
\multirow{3}{*}{20}    & Lower Bound        & 1058.3  & 1055.9 & 989.2  & 572.7  & 964.7  & 1051.8  \\
                       & Upper Bound        & 22563.0 & 5501.9 & 5487.4 & 700.4  & 1115.4 & 1212.6 \\
                       & Root Node Gap (\%) & 95.3    & 80.8   & 82.0   & 18.2   & 13.3   & 13.5   \\ \midrule
\multirow{3}{*}{23}    & Lower Bound        & 1260.9  & 1237.3 & 1285.3 & 1260.9 & 1219.0 & 1235.4 \\
                       & Upper Bound        & 32158.6 & 5683.3 & 5783.4 & 1495.7 & 1459.8 & 1436.2 \\
                       & Root Node Gap (\%) & 96.1    & 78.2   & 77.8   & 15.7   & 15.4   & 15.1   \\ \midrule
\multirow{3}{*}{25}    & Lower Bound        & 1326.2  & 1433.5 & 1388.2 & 675.1  & 1391.4 & 1389.9 \\
                       & Upper Bound        & 37700.7 & 5879.5 & 5886.3 & 854.3  & 1757.4 & 1753.6 \\
                       & Root Node Gap (\%) & 96.5    & 75.6   & 76.4   & 21.0   & 20.8   & 20.7   \\ \bottomrule
\end{tabular}
\end{table}

Comparing the results in Tables \ref{tab:Bounds_nw} and \ref{tab:Bounds_rc}, a few observations stand out. First, comparing TBD to BD, the upper bounds in TBD decrease substantially, with a reduction of 68.7\% compared to BD for the \textit{n\_w\_} and \textit{rc\_} instances.

Second, partial Benders decomposition helps tighten the lower bound at the root node, i.e., comparing BDP and TBDP with BD and TBD. It alleviates the primal inefficiencies with redundant solutions in early iterations of the BD and TBD methods. Adding any scenarios improves the lower and upper bound on average by 14.5\% and 86.1\%  over the \textit{n\_w\_} instances, respectively. For the \textit{rc\_} instances, the lower bound improves on average by 11.4\% with the addition of any scenarios. Similarly, the upper bound decreases by 81.8\% for the \textit{rc\_} instances.

We observe the benefit of choosing scenarios with { the} new scenario-retention strategy as it tightens the lower bound on average by 1.5\%  for \textit{n\_w\_} instances and 1.7\% for  \textit{rc\_} instances comparing the {BDS} and TBDS methods to the {BDP} and TBDP methods. Comparing TBDP with TBDS, we see that in TBDS, on average, the root node gap tightens by 3.4\%.

As the number of customers increases in both instances, the benefits gained from new scenario-retention strategy and the two-step decomposition become more pronounced, although the root gap remains high.
After examining the method's computational performance in depth, we shift our attention to problem-specific analysis. 

\section{Conclusions} \label{Sec:Conclusion}
{Assigning (accurate) time windows to customers in last-mile delivery can significantly improve customer satisfaction, but if done sub-optimally, it can increase costs tremendously. So far, the \textit{time-window assignment vehicle routing problem with stochastic travel times} (TWATSP-ST) has not been optimally solved. The TWATSP-ST is a fundamental and challenging combinatorial problem formulated as a two-stage stochastic mixed-integer program with continuous recourse for which no efficient exact solution methods exist.
We present a new method called Two-Step Benders Decomposition with Scenario Clustering (TBDS) for solving the TWATSP-ST, which can readily be applied to general two-stage stochastic mixed-integer programs.  }

Our method combines and generalizes the recent advancements in Benders decomposition and scenario clustering techniques. In particular, our TBDS method introduces a novel two-step decomposition strategy for the binary and continuous first-stage variables, resulting in improved continuous first-stage solutions while generating optimality cuts. This two-step decomposition approach leads to high-quality initial first-stage solutions, effectively reducing unnecessary iterations in current state-of-the-art Benders decomposition approaches. Consequently, it enhances computational efficiency by facilitating faster convergence.

The second key contribution of TBDS is the new scenario-retention strategy, which is, to the best of the authors' knowledge, the first time that such scenario clustering techniques and state-of-the-art Benders decomposition approaches are combined. By clustering the scenarios, we improve the linear programming (LP) relaxation of the master problem, obtaining superior lower bounds in the early iterations. By combining these essential elements (i.e., the two-step decomposition and the scenario clustering), our method achieves consistently tighter bounds at the root node and produces higher quality incumbent solutions compared to state-of-the-art approaches in the extant literature, including Benders dual decomposition and partial Benders decomposition. Specifically, these methods can be considered special cases of TBDS.

Extensive experimental results demonstrate the effectiveness of TBDS. Our method solves more instances to optimality, and significantly better lower and upper bounds are obtained for the instances not solved to optimality. In particular, TBDS achieves optimality for 73.9\% of the instances in our benchmark set, surpassing other benchmark algorithms that can solve at most 47\% of the instances. This showcases the superior performance and efficiency of our TBDS method in handling TWATSP-ST. Furthermore, our study reveals that the simultaneous optimization of time windows and routing leads to a noteworthy 11.1\% improvement in total costs while incurring only a minor increase in routing costs. Allowing different time window lengths enables hedging against high variances encountered throughout the route. As a result, our method produces shorter routes with fewer time window violations, contributing to the overall cost reduction.

In future studies, the alternating pattern of high and low variance arcs traveled, as observed in the structure of the optimal solution, can serve as a foundation for developing efficient heuristics explicitly tailored for the TWATSP-ST. Additionally, evaluating our algorithm on extended versions of the vehicle routing problem, such as the capacitated vehicle routing problem and the multi-depot vehicle routing problem, can offer valuable insights into the versatility and applicability of our TBDS method in diverse real-world scenarios.

From a methodological perspective, we envision that the problem-specific selection of representative scenarios combined with Benders decomposition approaches can start several new research lines. For instance, using supervised learning to predict which scenarios to label as representative based on instance-specific information such as vehicle information and customer locations seems promising. Mainly for applications with limited information, it would be valuable to research how well such predictions translate to slightly different settings. Alternatively, methods other than those in TBDS can be developed and tested to generate representative scenarios.

\ACKNOWLEDGMENT{
Albert H. Schrotenboer has received support from the Dutch Science Foundation (NWO) through grant VI.Veni.211E.043
}

\bibliographystyle{informs2014trsc} 
\bibliography{references.bib} 

\end{document}